\documentclass[11pt,letter]{article}

\usepackage[english]{babel}
\usepackage{amsmath}
\usepackage{amssymb}
\usepackage{bm}
\usepackage{mathrsfs}
\usepackage{graphics,graphicx,theorem}
\usepackage{dsfont}
\newtheorem{lemma}{Lemma}
\newtheorem{theorem}{Theorem}

{\theorembodyfont{\rmfamily}}
{\theorembodyfont{\rmfamily}}
{\theorembodyfont{\rmfamily}}

\allowdisplaybreaks  
\newcommand{\qed}{$\square$}

\makeatletter
\renewcommand\@biblabel[1]{}
\makeatother

\begin{document}
\begin{center}
\textmd{\LARGE{\bfseries{{Linear and quadratic functionals\\ of
random hazard rates: an asymptotic analysis}}}}
\end{center}

\medskip

\begin{center}
{\large Giovanni Peccati$^1$ and Igor Pr\"unster$^2$}
\end{center}

\begin{quote}
\begin{small}
\noindent $^1$ Laboratoire de Statistique Th\'eorique et
Appliqu\'ee, Universit\'e Paris VI, rue du Chevaleret 175, 75013
Paris, France. \\ \textit{e-mail:}
\texttt{giovanni.peccati@gmail.com}

\noindent $^2$ Dipartimento di Statistica e Matematica Applicata and
ICER, Universit\`a degli Studi di Torino, Piazza Arbarello 8, 10122
Torino, Italy. \\ \textit{e-mail}: \texttt{igor@econ.unito.it}

\end{small}
\end{quote}

\date{}

\begin{quote}
{\bf Abstract.} A popular Bayesian nonparametric approach to
survival analysis consists in modeling hazard rates as kernel
mixtures driven by a completely random measure. In this paper we
derive asymptotic results for linear and quadratic functionals of
such random hazard rates. In particular, we prove central limit
theorems for the cumulative hazard function and for the path--second
moment and path--variance of the hazard rate. Our techniques are
based on recently established criteria for the weak convergence of
single and double stochastic integrals with respect to Poisson
random measures. We illustrate our results by considering specific
models involving kernels and random measures commonly exploited in
practice.

\smallskip

{\bf Keywords}: Asymptotics; Bayesian Nonparametrics; Central limit
theorem; Path--variance; Random hazard rate; Survival analysis;
Completely random measure; Multiple Wiener-It\^{o} integral.

\smallskip

{\bf AMS 2000 classifications}: 62G20; 60G57.
\end{quote}

\section{Introduction}

Survival analysis has been the focus of many contributions to
Bayesian nonparametric theory and practice. Indeed, many statistical
problems arising in the framework of survival analysis require
function estimation and, hence, they are ideally suited for a
nonparametric treatment. Essentially, two closely related lines of
research have been pursued: the first is represented by the
introduction of models for the random cumulative distribution
function whereas the second deals with models for the random hazard
rate and the random cumulative hazard. As for the former most
proposals fall within the class of neutral to the right processes
due to Doksum (1974): among others, we mention Ferguson (1974),
Ferguson and Phadia (1979), Walker and Muliere (1997), Walker and
Damien (1998), Epifani, Lijoi and Pr\"{u}nster (2003), James (2006).
As for the latter, one can distinguish models leading to a
cumulative hazard which is almost surely discrete and models for
which it is almost surely absolutely continuous. The famous beta
process derived in Hjort (1990) belongs to the first class along the
contributions of, e.g., Kalbfleisch (1978), Kim (1999), Kim and Lee
(2003), De Blasi and Hjort (2006). The second class focuses on the
hazard rate which is modeled as a mixture and has recently received
much attention due to a relatively simple implementation in
applications. After the seminal papers of Dykstra and Laud (1981)
and Lo and Weng (1989), important developments dealing also with
more general multiplicative intensity models can be found in, Laud,
Smith and Damien (1996), Ibrahim, Chen and Mac Eachern (1999), James
(2003), Ishwaran and James (2004), Nieto--Barajas and Walker (2004,
2005), James (2005), Ho (2006), among others. Passing from a hazard
rate function to the corresponding model for the cumulative
distribution function is straightforward if the hazard rate is
almost surely absolutely continuous, but quite subtle otherwise. See
Hjort (1990) and James (2006), who establishes a nice link via the
notion of spatial neutral to the right process. It is also worth
noting that all models share a common feature, namely, that their
basic building block is represented by an increasing additive
process (see Sato, 1999) or more generally by a completely random
measure, a notion introduced in Kingman (1967).

Let us focus attention on hazard rates that are modeled as mixtures.
Denote by $U$ a positive absolutely continuous random variable
representing the lifetime and assume that its random hazard rate is
of the form
\begin{equation}
\tilde{h}(t)=\int_{\mathbb{X}}k(t,x)\tilde{\mu}(\mathrm{d}x),
\label{eq:mixture_hazard}
\end{equation}%
where $k$ is a kernel and $\tilde{\mu}$ a completely random measure
on some
space $\mathbb{X}$. The cumulative hazard is then given by $\tilde{H}%
(t)=\int_{0}^{t}\tilde{h}(s)\mathrm{d}s$. Note that, given $\tilde{\mu}$, $%
\tilde{h}$ represents the hazard rate of $U$, that is
\[
\tilde{h}(t)\ \mathrm{d}t=\mathds{P}(t\leq U\leq t+\mathrm{d}t|U\geq t,%
\tilde{\mu}).
\]%
From \eqref{eq:mixture_hazard}, provided $\tilde{H}(t)\rightarrow
\infty $ for $t\rightarrow \infty $ almost surely, one can define
a random density function $f$ as
\begin{equation}
\tilde{f}(t)=\tilde{h}(t)\ \exp (-\tilde{H}(t))  \nonumber
\end{equation}%
where $\tilde{S}(t):=\exp (-\tilde{H}(t))$ is the so--called
survival function providing the probability that $U> t$.
Consequently the random cumulative distribution function of $U$ is
of the form $\tilde{F}(t)=1-\exp (-\tilde{H}(t))$. Such models,
often referred to as life--testing models,
have been considered in Dykstra and Laud (1981) and Lo and Weng (1989) with $%
\tilde{\mu}$ being an extended gamma process, also known as weighted
gamma process. Nieto--Barajas and Walker (2004), instead, used a
weighted version
of a gamma compound Poisson process. Analysis beyond gamma--like choices of $%
\tilde{\mu}$ was not possible due to the lack of a suitable and
implementable posterior characterization: however, in James (2005)
this goal has been achieved and many choices for $\tilde{\mu}$ can
now be explored. See also Ho (2006) for a posterior characterization
via S--paths.

In this paper, we provide asymptotic results for random hazard rates
constructed via a mixture approach as in \eqref{eq:mixture_hazard}.
In particular, for $i=1,2,3$, we will be interested in establishing
the existence of two positive functions $\tau _{i}\left( T\right) $
and $\eta _{i}\left( T\right) $ such that the following Central
Limit Theorems (CLTs
in the sequel) take place as $T\rightarrow +\infty $:%
\begin{eqnarray}
&&\eta _{1}\left( T\right) \times \left[ \tilde{H}(T)-\tau
_{1}\left( T\right) \right] \overset{\mathrm{law}}{\longrightarrow
}X_{1}\left( \sigma
_{1}\right)   \label{CLT1} \\
&&\eta _{2}\left( T\right) \times \left[ \frac{1}{T}\int_{0}^{T}\tilde{h}%
(t)^{2}\mathrm{d}t-\tau _{2}\left( T\right) \right] \overset{\mathrm{law}}{%
\longrightarrow }X_{2}\left( \sigma _{2}\right)   \label{CLT2} \\
&&\eta _{3}\left( T\right) \times \left[
\frac{1}{T}\int_{0}^{T}\left[
\tilde{h}(t)-\tilde{H}(T)/T\right] ^{2}\mathrm{d}t-\tau _{3}\left( T\right) %
\right] \overset{\mathrm{law}}{\longrightarrow }X_{3}\left( \sigma
_{3}\right) \text{,}  \label{CLT3}
\end{eqnarray}%
where, for $i=1,2,3$, $X_{i}\left( \sigma _{i}\right) $ is a
centered Gaussian random variable, with variance $\sigma _{i}$
depending on the analytic structures of $\tilde{\mu}$ and $k$. For a
fixed $T>0$, the random
objects $T^{-1}\int_{0}^{T}\tilde{h}(t)^{2}\mathrm{d}t$ and $%
T^{-1}\int_{0}^{T}\left[ \tilde{h}(t)-\tilde{H}(T)/T\right]
^{2}\mathrm{d}t$ are called, respectively, the (realized)
\textsl{path-second moment} and the (realized)
\textsl{path-variance} associated with $\tilde{h}$. As we will
point out in the subsequent sections, weak convergence results such as (\ref%
{CLT1}), (\ref{CLT2}) and (\ref{CLT3}) give a description of the
overall variability of the hazard rate $\tilde{h}(t)$, by providing
a synthetic answer to the following crucial questions: (i)
\textquotedblleft How fast does the cumulative hazard rate diverge
from its long-term trend $\tau _{1}\left( T\right)
$?\textquotedblright , (ii) \textquotedblleft How fast increases the
magnitude of the fluctuations of $\tilde{h}(t)$ above
zero?\textquotedblright , and (iii) \textquotedblleft How big are
the oscillations of $\tilde{h}(t)$ around its average
value?\textquotedblright . To the authors knowledge, this represents
a completely new line of research. Indeed, by now, many results have
been obtained in terms of consistency of posterior distributions.
See Ghosh and Ramamoorthi (2003) for an exhaustive account. However,
little is known about the distributional behavior of the
prior ingredients of a Bayesian nonparametric model such as %
\eqref{eq:mixture_hazard}, in particular with reference to
functionals of statistical relevance. In the more conventional setup
of random probability measures, instead of the one concerning hazard
rates considered here, the first results on linear functionals of
the Dirichlet process were achieved in the pioneering paper of
Cifarelli and Regazzini (1990), whereas the variance functional is
studied in Cifarelli and Melilli (2000) and Regazzini, Guglielmi and
Di Nunno (2002). One may try to adopt the approach of Regazzini,
Lijoi and Pr\"{u}nster (2003) based on Gurland's inversion formula
to derive expressions for the distribution of linear functionals of
general random hazards as in \eqref{eq:mixture_hazard}, but to
tackle quadratic functionals seems impossible to date. In light of
these considerations, it seems important to remark that, despite the
theoretical relevance of our asymptotic results, they also turn out
to be helpful in terms of prior specification: on one hand they can
serve as a guide for deciding on which particular completely random
measure $\tilde{\mu}$ basing the model \eqref{eq:mixture_hazard} and
on the other hand, once $\tilde{\mu}$ is chosen, provide hints for
selecting the parameters of $\tilde{\mu}$. Indeed, up to now these
two steps were carried out in a conventional way, leaving aside the
problem of properly incorporating prior knowledge, in particular
with respect to the choice of $\tilde{\mu}$. A first contribution
highlighting the different clustering behaviors induced by
alternative random measures in the context of mixtures for Bayesian
density estimation is provided in Lijoi, Mena and Pr\"{u}nster
(2005). See also Ishwaran and James (2001).

The paper is structured as follows. In Section 2 we introduce some
basic concepts and notations. In Section 3 we state the main results
concerning linear and quadratic functionals of random hazard rates.
In particular, we derive CLTs for the cumulative hazard function and
for the path--second moment and path--variance of the hazard rate.
Moreover, we provide a useful comparison theorem which allows to
bypass the verification of the most delicate conditions thus leading
to obtain CLTs for hazard rates based on complex kernels or random
measures. Section 4 is devoted to applications: we consider specific
models involving kernels and random measures commonly exploited in
practice and analyze their asymptotic behavior in detail. In Section
5 the proofs of our results are provided and the techniques used to
establish them are illustrated. Section 5 contains some concluding
remarks together with possible extensions and an outline of future
work.

\section{Basic concepts and notations}

We start by introducing the main concepts and notations employed
throughout
the paper. Consider a measure space $(\mathbb{X},\mathscr{X})$, where $%
\mathbb{X}$ is a complete and separable metric space and
$\mathscr{X}$ is the usual Borel $\sigma $--field. Introduce a
\textsl{Poisson random measure}
$\tilde{N}$, defined on some probability space $(\Omega ,\mathscr{F},%
\mathds{P})$ and taking values in the set of non--negative counting
measures
on $(\mathds{R}^{+}\times \mathbb{X},\mathscr{B}(\mathds{R}^{+})\otimes %
\mathscr{X})$, with non-atomic \textsl{intensity measure} $\nu $,
i.e.
\begin{equation}
\mathds{E}\left[ \tilde{N}(\mathrm{d}v,\mathrm{d}x)\right] =\nu (\mathrm{d}v,%
\mathrm{d}x)  \nonumber
\end{equation}%
and, for any $A\in \mathscr{B}(\mathds{R}^{+})\otimes \mathscr{X}$
such that
$\nu (A)<\infty $, $\tilde{N}(A)$ is a Poisson random variable of parameter $%
\nu \left( A\right) $. Moreover, given any finite collection of
pairwise
disjoint sets, $A_{1},\ldots ,A_{k}$, in $\mathscr{B}(\mathds{R}^{+})\otimes %
\mathscr{X}$, the random variables $\tilde{N}(A_{1}),\ldots
,\tilde{N}(A_{k}) $ are mutually independent. Throughout the paper,
$\mathds{E}[\, \cdot \,]$ will denote expectation with respect to
$\mathds{P}$. Moreover, the intensity measure $\nu $ must satisfy
\begin{equation}
\int_{\mathds{R}^{+}}(v\wedge 1)\nu (\mathrm{d}v,\mathbb{X})<\infty
\nonumber
\end{equation}%
where $a\wedge b=min\{a,b\}$. See Daley and Vere--Jones (1988) for
an exhaustive account on Poisson random measures.

Recall that, according e.g. to Daley and Vere-Jones (1988), a Borel measure $%
\mu $ on some Polish space endowed with the Borel $\sigma $--algebra
is said to be \textsl{boundedly finite} if $\mu (A)<+\infty $ for
every bounded measurable set $A$. Let now
$(\mathbb{M},\mathscr{B}(\mathbb{M}))$ be the
space of boundedly finite measures on $(\mathbb{X},\mathscr{B}(\mathbb{X}))$%
. We suppose that $\mathbb{M}$ is equipped with the topology of weak
convergence and that $\mathscr{B}(\mathbb{M})$ is the corresponding Borel $%
\sigma $--field. Let $\tilde{\mu}$ be a random element, defined on $(\Omega ,%
\mathscr{F},\mathds{P})$ and with values in $(\mathbb{M},\mathscr{B}(\mathbb{%
M}))$, which can be represented as a linear functional of the
Poisson random measure $\tilde{N}$ (with intensity $\nu $) as
follows
\begin{equation}
\tilde{\mu}(B)=\int_{\mathds{R}^{+}\times B}s\,\tilde{N}(\mathrm{d}s,\mathrm{%
d}x)\qquad \forall B\in \mathscr{B}(\mathbb{X}).  \nonumber
\end{equation}%
It can be easily deduced from the properties of $\tilde{N}$ that
$\tilde{\mu} $ is, in the terminology of Kingman (1967), a
\textsl{completely random measure} (CRM) on $\mathbb{X}$, i.e.

\begin{itemize}
\item[(i)] $\tilde \mu(\emptyset)=0 \ $ a.s.-$\mathds{P}$

\item[(ii)] for any collection of disjoint sets in $\mathscr{B}(\mathbb{X})$%
, $B_1,B_2,\ldots$, the random variables $\tilde \mu(B_1),\tilde
\mu(B_2),\ldots$ are mutually independent and $\tilde \mu(\cup_{j\ge
1} B_j)=\sum_{j\ge 1}\tilde \mu(B_j)$ holds true a.s.-$\mathds{P}$.
\end{itemize}

Now let $\mathscr{G}_{\nu }$ be the space of functions $g:\mathbb{X}%
\rightarrow \mathds{R}^{+}$ such that $\int_{\mathds{R}^{+}\times \mathbb{X}%
}[1-\mathrm{e}^{-s\,g(x)}]$ $\nu (\mathrm{d}s,\mathrm{d}x)<\infty
$. Then, the law of $\tilde{\mu}$ is uniquely characterized by its
\textit{Laplace functional} which, for any $g$ in
$\mathscr{G}_{\nu }$, is given by
\begin{equation}
\mathds{E}\left[ \mathrm{e}^{-\int_{\mathbb{X}}g(x)\,\tilde{\mu}(\mathrm{d}%
x)}\right] =\exp \left\{ -\int_{\mathds{R}^{+}\times \mathbb{X}}[1-\mathrm{e}%
^{-s\,g(x)}]\ \nu (\mathrm{d}s,\mathrm{d}x)\right\}
\label{eq:Laplace}
\end{equation}%
For details and further references on CRMs see Kingman (1993). From %
\eqref{eq:Laplace} it is apparent that the law of the CRM
$\tilde{\mu}$ is completely determined by the corresponding
intensity measure $\nu $. This suggests a simple and useful
distinction of the random measures we deal with, according to the
decomposition of $\nu $. Letting $\lambda $ be a non--atomic and
$\sigma $--finite measure on $\mathbb{X}$, we have:

\begin{itemize}
\item[(a)] if $\nu(\mathrm{d} v,\mathrm{d} x)=\rho(\mathrm{d} v)\, \lambda(%
\mathrm{d} x)$, for some measure $\rho$ on $\mathds{R}^+$, we say
that the corresponding $\tilde N$ and $\tilde \mu$ are
\textit{homogeneous};

\item[(b)] if $\nu (\mathrm{d}v,\mathrm{d}x)=\rho (\mathrm{d}v|x)\,\lambda (%
\mathrm{d}x)$, where $\rho :\mathscr{B}(\mathds{R}^{+})\times \mathbb{X}%
\rightarrow \mathds{R}^{+}$ is a kernel (i.e. $x\mapsto \rho (C|x)$ is $%
\mathscr{B}(\mathbb{X})$--measurable for any $C\in \mathscr{B}(\mathds{R}%
^{+})$ and $\rho (\,\cdot \,|x)$ is a $\sigma $--finite measure on $%
\mathscr{B}(\mathds{R}^{+})$ for any $x$ in $\mathbb{X}$), we say
that the
corresponding random measures $\tilde{N}$ and $\tilde{\mu}$ are \textit{%
non--homogeneous}.
\end{itemize}

In the sequel we consider CRM $\tilde \mu$ whose intensity measures
satisfy
\begin{equation}  \label{eq:H1}
\int_{\mathds{R}^{+}\times \mathbb{X}}\rho (dv|x)\,
\lambda(dx)=+\infty \tag{H1}
\end{equation}

In the homogeneous case, \eqref{eq:H1} reduces to $\max \{\,\rho (\mathds{R}%
^{+})\,;\,\lambda (\mathbb{X})\,\}=+\infty $, which is tantamount of
requiring either infinite activity of $\tilde{\mu}$ i.e.
$\tilde{\mu}$
jumping infinitely often on any bounded $A\in \mathscr{X}$ or to consider $%
\tilde{\mu}$ with unbounded support $S$ such that $\lambda \left(
S\right) =+\infty .$ In the non--homogeneous case, for \eqref{eq:H1}
to hold it is enough that $\tilde{\mu}$ jumps infinitely often on
some bounded set of positive $\lambda $--measure. It is clear that
\eqref{eq:H1} is met by the CRM usually considered in the
literature. In the subsequent sections, as illustrations of our
general results, we will consider the following CRMs:

\begin{enumerate}
\item Generalized gamma CRM: its intensity measure is homogeneous and given
by
\begin{equation}
\nu (\mathrm{d}v,\mathrm{d}x)=\frac{1}{\Gamma (1-\sigma )}\,\frac{\mathrm{e}%
^{-\gamma v}}{v^{1+\sigma}}\, \mathrm{d}v\, \lambda (\mathrm{d}x)
\label{eq:gg}
\end{equation}%
where $\sigma \in (0,1)$ and $\gamma >0$. This class, studied in
Brix (1999), can be characterized as the tilted exponential family
generated by the positive stable laws. It includes the inverse
Gaussian CRM for $\sigma =1/2$ and the gamma CRM as $\sigma
\rightarrow 0$.

\item Extended gamma CRM: its non--homogeneous intensity measure is of the
form
\begin{equation}
\nu (\mathrm{d}v,\mathrm{d}x)=\frac{\mathrm{e}^{-\beta (x)v}}{v}\, \mathrm{d}%
v\, \lambda (\mathrm{d}x)  \label{eq:extended_gamma}
\end{equation}%
where $\beta $ is a strictly positive function on $\mathbb{X}$. This
class dates back to Dykstra and Laud (1981) and Lo and Weng (1989).
The gamma CRM arises if $\beta $ is constant.

\item Beta CRM: its non--homogeneous intensity measure is given by
\begin{equation}
\nu (\mathrm{d}v,\mathrm{d}x)=\mathbb{I}_{(0,1)}(v)\,c(x)\,\frac{%
(1-v)^{c(x)-1}}{v}\ \mathrm{d}v\,\lambda (\mathrm{d}x)
\label{eq:beta}
\end{equation}%
where $c$ is some strictly positive function on $\mathbb{X}$ and $\mathbb{I}%
_{A}$ stands for the indicator function of set $A$. Note that the
class of beta CRM, which is due to Hjort (1990), has the
particularity of allowing only jumps of sizes less than 1.
\end{enumerate}

Having settled the basics regarding the background driving CRM in %
\eqref{eq:mixture_hazard}, we now have to define the kernel: $k$ is
a
jointly measurable application from $\mathds{R}^{+}\times \mathbb{X}$ to $%
\mathds{R}^{+}$, such that $\int_{\mathbb{X}}k(t,x)\lambda (\mathrm{d}%
x)<+\infty $ and $\int_{\,\cdot \,}k(t|x)\mathrm{d}t$ is a $\sigma
$--finite measure on $\mathscr{B}(\mathds{R}^{+})$ for any $x$ in
$\mathbb{X}$. Given these two ingredients the random hazard rate in
\eqref{eq:mixture_hazard} is properly defined.

A further technical assumption we will make throughout the paper is
represented by the following conditions
\begin{align}
& {\int_{\mathds{R}^{+}\times \mathbb{X}}k}\left( t,x\right) ^{j}v^{j}\ {%
\rho (\mathrm{d}v|x)\,\lambda (\mathrm{d}x)<+\infty }\quad \forall
t,\text{
\ }j=1,2,4;  \tag{H2}  \label{eq:H2} \\
& \int_{\left[ 0,T\right] }{\int_{\mathds{R}^{+}\times
\mathbb{X}}k}\left(
t,x\right) ^{j}v^{j}\ {\rho (\mathrm{d}v|x)\,\lambda (\mathrm{d}x)\,\mathrm{d%
}t<+\infty }\quad \forall T>0,\text{ \ }j=1,2,4.  \nonumber
\end{align}%
If, for $j=1,2,4$, the application $x\mapsto \int_{\mathds{R}^{+}}v^{j}\rho (%
\mathrm{d}v|x)$ is\textsl{\ bounded} by some finite constant (which
is typically the case), then the first condition in \eqref{eq:H2}
reduces to requiring that the function $x\mapsto {k}\left(
t,x\right) ^{j}$ is integrable with respect to $\lambda $ for every
$t$, whereas the second line of \eqref{eq:H2} boils down to the
assumption that the application $\left( t,x\right) \mapsto k\left(
t,x\right) $ is an element of $\cap _{j=1,2,4}L^{j}(\left[
0,T\right] \times \mathbb{X},\,\mathrm{d}t\,\lambda \left(
\mathrm{d}x\right) )$ for every $T>0$. Hence, in the uniformly
bounded case \eqref{eq:H2} is a condition not involving the CRM, but
just
the kernel. Moreover, it is easy to see that the quantity $\int_{\mathds{R}%
^{+}}v^{j}\rho (\mathrm{d}v|x)$, $j=1,2,4$, is bounded in $x$
whenever $\rho (\mathrm{d}v|x)$ is associated to one of the three
classes of CRMs defined above (see (\ref{eq:gg}),
(\ref{eq:extended_gamma}) and (\ref{eq:beta})). We
shall also note that, in the homogeneous case, \eqref{eq:H2} implies that $%
\int_{\mathds{R}^{+}}v^{j}\rho (\mathrm{d}v)<+\infty $, $j=1,2,4$.
An example of a homogeneous CRM, which does not meet \eqref{eq:H2}
is the
stable CRM for which $\rho (\mathrm{d}v)=v^{-1-\sigma }\mathrm{d}v$ and $%
\sigma \in (0,1)$. Note that the stable CRM can be recovered from
the generalized gamma class by allowing $\gamma =0$ in
\eqref{eq:gg}: we have excluded this possibility since it does not
meet \eqref{eq:H2}.

\subsection{Further notation}

For $q,p\geq 1$, we note
\[
L^{p}\left( \nu ^{q}\right) =L^{p}((\mathds{R}^{+}\times \mathbb{X)}^{q},(%
\mathscr{B}(\mathds{R}^{+})\otimes \mathscr{X})^{q},\nu ^{q})
\]%
the Banach space of real-valued functions $f$ on
$(\mathds{R}^{+}\times \mathbb{X)}^{q}$, such that $\left\vert
f\right\vert ^{p}$ is integrable with respect to $\nu ^{q}$ $:=$
$\nu ^{\otimes q}$. We will systematically write $L^{p}\left( \nu
^{1}\right) =L^{p}\left( \nu \right) $ for $p\geq 1$. The symbol
$L_{s}^{2}\left( \nu ^{2}\right) $ is used to denote the subspace
of $L^{2}\left( \nu ^{2}\right) $ composed of \textsl{symmetric functions }on%
\textsl{\ }$(\mathds{R}^{+}\times \mathbb{X)}^{2}$. By symmetric, we
mean that every $f\in L_{s}^{2}\left( \nu ^{2}\right) $ is such that
$f\left( s,x;t,y\right) =f\left( t,y;s,x\right) $ for every $\left(
s,x\right)
,\left( t,y\right) \in \mathds{R}^{+}\times \mathbb{X}.$ We also write $%
L_{s,0}^{2}\left( \nu ^{2}\right) $ to indicate the subset of $%
L_{s}^{2}\left( \nu ^{2}\right) $ composed of symmetric functions \textsl{%
vanishing on diagonals}, i.e. such that their support is contained
in the \textsl{purely non-diagonal set} $D_{0}^{2}=\{\left(
s,x;t,y\right) :\left( s,x\right) \neq \left( t,y\right) \}$.

We now turn to the definition of three basic auxiliary kernels which
are
associated to a given $f\in L_{s}\left( \nu ^{2}\right) $: (i) the kernel $%
f\star _{1}^{0}f$ is defined on $(\mathds{R}^{+}\times
\mathbb{X)}^{3}$ and
is given by%
\begin{equation}
f\star _{1}^{0}f\left( t_{1},x_{1};t_{2},x_{2};t_{3},x_{3}\right)
=f\left(
t_{1},x_{1};t_{2},x_{2}\right) f\left( t_{2},x_{2};t_{3},x_{3}\right) \text{%
; }  \label{eq : f10}
\end{equation}%
(ii) $f\star _{1}^{1}f$ is defined on $(\mathds{R}^{+}\times
\mathbb{X)}^{2}$
and is actually a \textsl{contraction} equal to%
\begin{equation}
f\star _{1}^{1}f\left( t_{1},x_{1};t_{2},x_{2}\right) =\int_{\mathds{R}%
^{+}\times \mathbb{X}}f\left( t_{1},x_{1};s,y\right) f\left(
s,y;t_{2},x_{2}\right) \nu \left( \mathrm{d}s,\mathrm{d}y\right) ;
\label{eq :f11}
\end{equation}%
(iii) $f\star _{2}^{1}f$ is defined on $(\mathds{R}^{+}\times
\mathbb{X)}$
and is given by%
\begin{equation}
f\star _{2}^{1}f\left( t,x\right) =\int_{\mathds{R}^{+}\times \mathbb{X}%
}f\left( t,x;s,y\right) ^{2}\nu \left(
\mathrm{d}s,\mathrm{d}y\right) . \label{eq : f21}
\end{equation}

\noindent Note that, by the Cauchy-Schwarz inequality and by the
symmetry and square-integrability of $f$, the kernel $f\star
_{1}^{1}f$ is necessarily an element of $L_{s}^{2}\left( \nu
^{2}\right) $. The three kernels defined above are the fundamental
building blocks to obtain explicit expressions for the moments and
the cumulants of the linear and quadratic functionals associated
with random hazard rates (when they exist). Such expressions enter
implicitly in the statements of the subsequent results, and are
mainly of a combinatorial nature. We refer the reader to Rota and
Wallstrom (1997) for an exhaustive analysis of the combinatorial
machinery underlying the construction of stochastic integrals with
respect to completely random measures.

In the subsequent sections it will be often convenient to work with
the \textsl{compensated} \textsl{Poisson random measure} canonically
associated
to $\tilde{N}$. Such an object is indicated by%
\begin{equation}
\tilde{N}^{c}=\left\{ \tilde{N}^{c}\left( A\right) :A\in \mathscr{B}(%
\mathds{R}^{+})\otimes \mathscr{X}\right\} ,  \label{eq : Comp P}
\end{equation}%
and is formally defined as the unique CRM on $(\mathds{R}^{+}\times \mathbb{X%
},\mathscr{B}(\mathds{R}^{+})\otimes \mathscr{X})$ such that
\begin{equation}
\tilde{N}^{c}\left( A\right) =\tilde{N}\left( A\right) -\nu \left(
A\right) \label{eq : Comp P2}
\end{equation}%
for every set $A$ of finite $\nu $-measure. For every $g\in
L^{2}\left( \nu
\right) $, we denote by%
\[
\tilde{N}^{c}\left( g\right) =\int_{\mathds{R}^{+}\times
\mathbb{X}}g\left( s,x\right) \tilde{N}^{c}\left(
\mathrm{d}s,\mathrm{d}x\right)
\]%
the Wiener-It\^{o} integral of $g$ with respect to $\tilde{N}^{c}$.
We recall that, for every $g\in L^{2}\left( \nu \right) $,
$\tilde{N}^{c}\left( g\right) $ is a centered and square integrable
random variable with an
infinitely divisible law, such that, for every $\lambda \in \mathbb{R}$,%
\begin{equation}
\mathds{E}\left[ \ \mathrm{e}^{i\lambda \tilde{N}^{c}\left( g\right) }\,%
\right] =\exp \left\{ \int_{\mathds{R}^{+}\times \mathbb{X}}\left[ \,\mathrm{%
e}^{i\lambda g\left( s,x\right) }-1-i\lambda g\left( s,x\right)
\right] \nu \left( \mathrm{d}s,\mathrm{d}x\right) \right\}
\label{PLK}
\end{equation}%
(compare with (\ref{eq:Laplace})). Moreover, for every $f,g\in
L^{2}\left(
\nu \right) $, one has the \textsl{isometric property}%
\begin{equation}
\mathds{E}\left[ \,\tilde{N}^{c}\left( f\right) \tilde{N}^{c}\left( g\right) %
\right] =\int_{\mathds{R}^{+}\times \mathbb{X}}f\left( s,x\right)
g\left( s,x\right) \nu \left( \mathrm{d}s,\mathrm{d}x\right)
:=\left( f,g\right) _{L^{2}\left( \nu \right) }.  \label{eq :
isoPoiss}
\end{equation}%
Note that (\ref{eq:Laplace}), (\ref{PLK}) and the isometric property (\ref%
{eq : isoPoiss}) imply that, for every $g\in L^{2}\left( \nu \right)
\cap
L^{1}\left( \nu \right) $,%
\begin{eqnarray}
\mathds{E}\left[ \,\tilde{N}\left( g\right) \right] &=&\int_{\mathds{R}%
^{+}\times \mathbb{X}}g\left( s,x\right) \nu \left( \mathrm{d}s,\mathrm{d}%
x\right) \text{ }  \label{eq : EN} \\[5pt]
\mathrm{Var}\!\left[ \,\tilde{N}\left( g\right) \right] &=&\mathrm{Var}\!%
\left[ \,\tilde{N}^{c}\left( g\right) \right]
=\int_{\mathds{R}^{+}\times
\mathbb{X}}g\left( s,x\right) ^{2}\nu \left( \mathrm{d}s,\mathrm{d}x\right) .%
\text{ }  \label{eq : VN}
\end{eqnarray}

\smallskip

\section{Main results: CLTs for linear and quadratic functionals}

In what follows, we shall develop several techniques, allowing to
study the asymptotic behavior of linear and quadratic functionals
associated to the random hazard rate $\tilde{h}(t)$ appearing in
(\ref{eq:mixture_hazard}). Concerning quadratic functionals, we will
be mainly interested in the path-variance and the path second moment
of $\tilde{h}(t)$. As will be clarified in Section \ref{S : Proofs},
our approach exploits the fact that any quadratic functional of
$\tilde{h}$ can be (uniquely) represented as a linear combination of
its expectation and of the following two random elements: (i)
the stochastic integral of a deterministic kernel with respect to $\tilde{N}%
^{c}$, and (ii) the double Wiener-It\^{o} integral of a
deterministic bivariate kernel with respect to the stochastic
product measure associated to $\tilde{N}^{c}$. According to the
results proved in Peccati and Taqqu (2006b) (see Section \ref{SS :
2ble int CLT}), the joint (weak) convergence of single and double
Poisson integrals can be characterized in terms of the asymptotic
negligibility of deterministic contraction kernels. We will show
that such contractions are indeed explicit functionals of the kernel
$k$ defining $\tilde{h}$. We shall first state the main general
results of the paper, and then describe in detail several
applications. The proofs are deferred to Section \ref{S : Proofs}.

\smallskip

Consider the random hazard rate $\tilde{h}$ defined in formula (\ref%
{eq:mixture_hazard}), and assume that the intensity of the
underlying Poisson CRM $\tilde{N}$ verifies (\ref{eq:H1}), and that
the positive kernel $k$ satisfies (\ref{eq:H2}). Moreover, for every
$T>0$ define the kernel
\begin{equation}
k_{T}^{\left( 0\right) }\left( s,x\right) =s\int_{0}^{T}k\left(
t,x\right) \mathrm{d}t,\text{ \ \ }\left( s,x\right) \in
\mathds{R}^{+}\times \mathbb{X} .  \label{kaT1}
\end{equation}%
Our first result concerns the asymptotic behavior of the cumulative
hazard rate $\tilde{H}(T)=\int_{0}^{T}\tilde{h}(t)\mathrm{d}t$.

\begin{theorem}
\label{T : cumHRasy} Suppose that: (i) $k_{T}^{\left( 0\right) }\in
L^{3}\left( \nu \right) $ for every $T$, and (ii) there exists a
strictly
positive function $T\mapsto C_{0}\left( k,T\right) $, such that, as $%
T\rightarrow +\infty $,
\begin{eqnarray}
C_{0}^{2}\left( k,T\right) \times \int_{\mathds{R}^{+}\times \mathbb{X}}%
\left[ k_{T}^{\left( 0\right) }\left( s,x\right) \right] ^{2}\nu
\left( \mathrm{d}s,\mathrm{d}x\right) &\rightarrow &\sigma
_{0}^{2}\left( k\right)
>0\text{, }  \label{a-} \\
C_{0}^{3}\left( k,T\right) \times \int_{\mathds{R}^{+}\times \mathbb{X}}%
\left[ k_{T}^{\left( 0\right) }\left( s,x\right) \right] ^{3}\nu
\left( \mathrm{d}s,\mathrm{d}x\right) &\rightarrow &0.  \label{aa}
\end{eqnarray}%
Then,
\begin{equation}
C_{0}\left( k,T\right) \times \left[ \tilde{H}(T)-\mathds{E}[\tilde H(T)]%
\right] \overset{\mathrm{law}}{\longrightarrow}X,  \label{aaa}
\end{equation}%
where $X\sim \mathscr{N}\left( 0,\sigma _{0}^{2}\left( k\right)
\right) $
\end{theorem}

Note that conditions (\ref{a-})-(\ref{aa}) only involve the analytic
form of the kernel $k$, and do not make any use the of the
asymptotic properties of the law of the process $\tilde{h}(t)$, such
as e.g. mixing. We now focus on the limiting behavior of the
quadratic functionals associated to the random hazard rate
$\tilde{h}$. To this end, we associate to $k\left( \cdot ,\cdot
\right) $, and to each $T>0$, the three auxiliary kernels:%
\begin{eqnarray}
k_{T}^{\left( 1\right) }\left( s,x;t,y\right) &=&\frac{st}{T}%
\int_{0}^{T}k\left( u,x\right) k\left( u,y\right) \mathrm{d}u,  \label{k1} \\
k_{T}^{\left( 2\right) }\left( s,x\right) &=&\frac{s^{2}}{T}%
\int_{0}^{T}k\left( u,x\right) ^{2}\mathrm{d}u,  \label{k2} \\
k_{T}^{\left( 3\right) }\left( s,x\right)
&=&\int_{\mathds{R}^{+}\times
\mathbb{X}}k_{T}^{\left( 1\right) }\left( s,x;u,w\right) \nu \left( \mathrm{d%
}u,\mathrm{d}w\right) .  \label{k3}
\end{eqnarray}
The kernel $k_{T}^{\left( 2\right) }$ can be obtained by restricting $%
k_{T}^{\left( 1\right) }$ to the \textsl{diagonal set}
$\{(s,x;t,y):\left( s,x\right) =(t,y)\}$. We will see in Section
\ref{S : Proofs} that the kernels $k_{T}^{\left( \cdot \right) }$
are intimately related to the objects defined in formulae (\ref{eq :
f10})-(\ref{eq : f21}). Note that, due to assumption (\ref{eq:H2})
and the Jensen and Cauchy-Schwarz inequalities, $k_{T}^{\left(
1\right) }\in L_{s}^{2}\left( \nu ^{2}\right) \cap L^{4}\left( \nu
^{2}\right) $, and also $k_{T}^{\left( 2\right) }\in L^{2}\left( \nu
\right) $. The following theorem provides a CLT for the path--second
moment of random hazard rates.

\begin{theorem}
\label{T : path2ndM} Suppose that $k_{T}^{\left( 3\right) }\in
L^{2}\left( \nu \right) \cap L^{1}\left( \nu \right) $,
$k_{T}^{\left( 2\right) }\in L^{3}\left( \nu \right) $ and that
there exists a strictly positive function $C_{1}\left( k,T\right) $
such that the following asymptotic conditions are satisfied as
$T\rightarrow +\infty $:

\begin{enumerate}
\item[$1.$] \ \ $2C_{1}^{2}\left( k,T\right) \left\Vert k_{T}^{\left(
1\right) }\right\Vert _{L^{2}\left( \nu ^{2}\right) }^{2}\rightarrow
\sigma _{1}^{2}\left( k\right) >0;$

\item[$2.$] \ \ $C_{1}^{4}\left( k,T\right) \left\Vert k_{T}^{\left(
1\right) }\right\Vert _{L^{4}\left( \nu ^{2}\right) }^{4}\rightarrow
0;$

\item[$3.$] \ \ $C_{1}^{4}\left( k,T\right) \left\Vert k_{T}^{\left(
1\right) }\star _{1}^{1}k_{T}^{\left( 1\right) }\right\Vert
_{L^{2}\left( \nu ^{2}\right) }^{2}\rightarrow 0;$

\item[$4.$] \ \ $C_{1}^{4}\left( k,T\right) \left\Vert k_{T}^{\left(
1\right) }\star _{2}^{1}k_{T}^{\left( 1\right) }\right\Vert
_{L^{2}\left( \nu \right) }^{2}\rightarrow 0;$

\item[$5.$] \ \ $C_{1}^{2}\left( k,T\right) \left\Vert k_{T}^{\left(
2\right) }+2k_{T}^{\left( 3\right) }\right\Vert _{L^{2}\left( \nu
\right) }^{2}\rightarrow \sigma _{2}^{2}\left( k\right) >0;$

\item[$6.$] \ \ $C_{1}^{3}\left( k,T\right) \left\Vert k_{T}^{\left(
2\right) }+2k_{T}^{\left( 3\right) }\right\Vert _{L^{3}\left( \nu
\right) }^{3}\rightarrow 0.$
\end{enumerate}

Then,
\begin{equation}
C_{1}\left( k,T\right) \times \left\{ \frac{1}{T}\int_{0}^{T}\tilde{h}(t)^{2}%
\mathrm{d}t-\frac{1}{T}\int_{0}^{T}\mathds{E}\lbrack \tilde{h}(t)^{2}]%
\mathrm{d}t\right\} \overset{\mathrm{law}}{\longrightarrow} X
\label{eq : 2ndMclt}
\end{equation}%
where $X\sim \mathscr{N}\left( 0,\sigma _{1}^{2}\left( k\right)
+\sigma _{2}^{2}\left( k\right) \right) $.
\end{theorem}

Note that
\begin{eqnarray*}
\left\Vert k_{T}^{\left( 3\right) }\right\Vert _{L^{1}\left( \nu
\right) }
&=&\int_{\mathds{R}^{+}\times \mathbb{X}}\int_{\mathds{R}^{+}\times \mathbb{X%
}}k_{T}^{\left( 1\right) }\left( s,x;u,w\right) \nu \left( \mathrm{d}u,%
\mathrm{d}w\right) \nu \left( \mathrm{d}s,\mathrm{d}x\right) \\
&=&\frac{1}{T}\int_{0}^{T}\left( \int_{\mathds{R}^{+}\times \mathbb{X}%
}sk\left( t,x\right) \nu \left( \mathrm{d}s,\mathrm{d}x\right) \right) ^{2}%
\mathrm{d}t.
\end{eqnarray*}%
Also, by applying formulae (\ref{eq : EN}) and (\ref{eq : VN}) (for every $%
t>0$) in the case $h\left( s,x\right) =sk\left( t,x\right) $, one
obtains
that%
\begin{eqnarray}
\frac{1}{T}\int_{0}^{T}\mathds{E}\lbrack \tilde{h}(t)^{2}]\mathrm{d}t &=&%
\frac{1}{T}\int_{0}^{T}\left( \int_{\mathds{R}^{+}\times
\mathbb{X}}sk\left( t,x\right) \nu \left(
\mathrm{d}s,\mathrm{d}x\right) \right) ^{2}\mathrm{d}t
\label{eq : semplVar} \\
&&\qquad +\frac{1}{T}\int_{0}^{T}\int_{\mathds{R}^{+}\times \mathbb{X}%
}s^{2}k\left( t,x\right) ^{2}\nu \left(
\mathrm{d}s,\mathrm{d}x\right) \mathrm{d}t.  \nonumber
\end{eqnarray}

The next theorem combines Theorem \ref{T : cumHRasy}\ and Theorem
\ref{T : path2ndM} to deal with path-variances of random hazard
rates.

\begin{theorem}
\label{P : path variance}Suppose that $\tilde{h}$ is such that
assumptions
\textnormal{\eqref{a-}--\eqref{aa}} are verified, and conditions \textnormal{%
1.--6.} of Theorem \textnormal{\ref{T : path2ndM}} are satisfied. If
there
exists a constant $\delta \left( k\right) \geq 0$ such that, as $%
T\rightarrow +\infty $,

\begin{enumerate}
\item[$1.$] \ \ $C_{1}\left( k,T\right) /\left( TC_{0}\left( k,T\right)
\right) ^{2}\rightarrow 0;$

\item[$2.$] \ \ $2C_{1}\left( k,T\right) \mathds{E}\lbrack \tilde{H}%
(T)]/\left( T^{2}C_{0}\left( k,T\right) \right) \rightarrow \delta
\left( k\right) ;$

\item[$3.$] \ \ $\left\Vert \,C_{1}\left( k,T\right) \left( k_{T}^{\left(
2\right) }+2k_{T}^{\left( 3\right) }\right) -\delta \left( k\right)
C_{0}\left( k,T\right) k_{T}^{\left( 0\right) }\right\Vert
_{L^{2}\left( \nu \right) }^{2}\rightarrow \sigma _{3}^{2}\left(
k\right) \geq 0,$
\end{enumerate}

then,
\begin{align}
& C_{1}\left( k,T\right) \!\times \!\left\{ \!\frac{1}{T}\!\int_{0}^{T}\!%
\left[ \tilde{h}(t)-\frac{\tilde{H}(T)}{T}\right] ^{2}\!\!\mathrm{d}t\!-\!%
\frac{1}{T}\!\int_{0}^{T}\!\!\mathds{E}\lbrack \tilde{h}(t)^{2}]\mathrm{d}t+%
\frac{\mathds{E}\lbrack \tilde{H}(T)]^{2}}{T^{2}}\!\right\}
\label{eq PVarclt} \\
=\,& C_{1}\left( k,T\right) \!\times \!\left\{ \!\frac{1}{T}\!\int_{0}^{T}\!%
\left[ \tilde{h}(t)-\frac{\tilde{H}(T)}{T}\right] ^{2}\!\!\mathrm{d}t\!-\!%
\frac{1}{T}\!\int_{0}^{T}\!\!\mathds{E}\left[ \tilde{h}(t)-\frac{\mathds{E}(%
\tilde{H}(T))}{T}\right] ^{2}\!\!\mathrm{d}t\!\right\}   \nonumber \\
\overset{\mathrm{law}}{\longrightarrow }& \ X,  \nonumber
\end{align}%
where $X\sim \mathscr{N}\left( 0,\sigma _{1}^{2}\left( k\right)
+\sigma _{3}^{2}\left( k\right) \right) $.
\end{theorem}

In view of (\ref{eq : VN}), one also has that%
\begin{equation}
\frac{1}{T}\int_{0}^{T}\mathrm{Var}\left( \tilde{h}(t)\right) \mathrm{d}t=%
\frac{1}{T}\int_{0}^{T}\int_{\mathds{R}^{+}\times
\mathbb{X}}s^{2}k\left( t,x\right) ^{2}\nu \left(
\mathrm{d}s,\mathrm{d}x\right) \mathrm{d}t. \nonumber
\end{equation}%
To conclude this subsection, we state a useful \textsl{comparison theorem }%
for random hazard rates. To this end, consider two completely random
Poisson
measures (on $\mathds{R}^{+}\times \mathbb{X}$) $\overline{N}$ and $%
\overline{\overline{N}} $, as well as positive kernels $\overline{k}$ and $%
\overline{\overline{k}}$. The $\sigma $-finite and non-atomic
intensity
measures of $\overline{N}$ and $\overline{\overline{N}}$ are denoted by $%
\overline{\nu }$ and $\overline{\overline{\nu }}$, respectively. We
assume
that $\overline{\nu }$ and $\overline{\overline{\nu }}$ both verify (\ref%
{eq:H1}), and that $\overline{k}$ and $\overline{\overline{k}}$ satisfy (\ref%
{eq:H2}). Finally, we suppose that, for every $B\in (\mathscr{B}(\mathds{R}%
^{+})\otimes \mathscr{X})$,%
\begin{equation}
\overline{\nu }\left( B\right) \leq \overline{\overline{\nu }}\left(
B\right) ,  \nonumber
\end{equation}%
and, for every $\left( t,x\right) \in \mathds{R}^{+}\times
\mathbb{X}$,
\begin{equation}
\overline{k}\left( t,x\right) \leq \overline{\overline{k}}\left(
t,x\right) \text{.}  \nonumber
\end{equation}
Throughout the paper, for strictly positive sequences $\left\{
a_{n}\right\}
$ and $\left\{ b_{n}\right\} $, we write $a_{n}\sim b_{n}$ if there exists $%
c\in \left( 0,+\infty \right) $ such that $a_{n}/b_{n}\rightarrow c$, as $%
n\rightarrow \infty$.

\begin{theorem}
\label{T : Comparison}Suppose that the pair $\left( \nu ,k\right) $
entering
the definition of the random hazard $\tilde{h}$ in \textnormal{%
\eqref{eq:mixture_hazard}} is such that, for every $B\in (\mathscr{B}(%
\mathds{R}^{+})\otimes \mathscr{X})$, $\overline{\nu }\left( B\right) $ $%
\leq $ $\nu \left( B\right) $ $\leq $ $\overline{\overline{\nu
}}\left( B\right) $ and, for every $\left( t,x\right) \in
\mathds{R}^{+}\times
\mathbb{X}$, $\overline{k}\left( t,x\right) $ $\leq $ $k\left( t,x\right) $ $%
\leq $ $\overline{\overline{k}}\left( t,x\right) $. Then, the
following three comparison criteria hold.

\textnormal{(A)} Assume that the two kernels $\overline{k}$ and $\overline{%
\overline{k}}$, with $\overline{\nu }$ and $\overline{\overline{\nu
}}$
substituting $\nu$, satisfy the conditions \textnormal{\eqref{a-}--\eqref{aa}%
} for some appropriate positive functions $C_{0}(\overline{k},T)$ and $C_{0}(%
\overline{\overline{k}},T)$ and constants $\sigma
_{0}^{2}(\overline{k})$
and $\sigma _{0}^{2}(\overline{\overline{k}})$. Suppose also that $C_{0}(%
\overline{k},T)$ $\sim $ $C_{0}(\overline{\overline{k}},T)$, and
consider a
positive function $C_{0}(k,T)$ such that $C_{0}(k,T)$ $\sim $ $C_{0}(%
\overline{k},T)$. Then, for every diverging sequence
$T_{n}\rightarrow +\infty $, there exists a subsequence
$T_{n^{\prime }}$ such that the CLT
\textnormal{\eqref{aaa}} holds as $n^{\prime }\rightarrow +\infty $, with $%
T_{n^{\prime }}$ substituting $T$, where $X$ is a centered Gaussian
random variable whose variance depends on the choice of $C_{0}\left(
k,T\right) $ and on $n^{\prime }$.

\textnormal{(B)} Assume that $\overline{k}$ and
$\overline{\overline{k}}$, with $\overline{\nu }$ and
$\overline{\overline{\nu }}$ substituting $\nu $, satisfy conditions
\textnormal{1.--6.} of Theorem \textnormal{\ref{T :
path2ndM}} for some positive functions $C_{1}(\overline{k},T)\ $and $C_{1}(%
\overline{\overline{k}},T)$ and constants $\sigma
_{j}^{2}(\overline{k})$ and $\sigma
_{j}^{2}(\overline{\overline{k}})$, $j=1,2$. Assume, moreover, that
$C_{1}(\overline{k},T)$ $\sim $ $C_{1}(\overline{\overline{k}},T)$,
and select a positive function $C_{1}\left( k,T\right) $ such that
$C_{1}\left(
k,T\right) $ $\sim $ $C_{1}(\overline{k},T)$. Then, for every sequence $%
T_{n}\rightarrow +\infty $, there exists a subsequence $T_{n^{\prime
}}$
such that the CLT \textnormal{\eqref{eq : 2ndMclt}} is verified (for $%
n^{\prime }\rightarrow +\infty $ and with $T_{n^{\prime }}$ substituting $T$%
) where $X$ is a centered Gaussian random variable whose variance
depends on $C_{1}\left( k,T\right) $ and $n^{\prime }$.

\textnormal{(C)} Suppose that $\overline{k}$, $\overline{\overline{k}}$, $%
C_{j}(\overline{k},T),$ $C_{j}(\overline{\overline{k}},T)$ and $C_{j}(k,T)\ $%
($j=0,1$) satisfy the assumptions pinpointed in Parts
\textnormal{(A)} and
\textnormal{(B)}, and suppose that they also meet the Conditions \textnormal{%
1.--3.} of Theorem \textnormal{\ref{P : path variance}}. Then, for
every
sequence $T_{n}\rightarrow +\infty $, there exists a subsequence $%
T_{n^{\prime }}$ such that the CLT \textnormal{\eqref{eq PVarclt}}
holds, for $n^{\prime }\rightarrow +\infty $ and with $T_{n^{\prime
}}$ substituting $T$.
\end{theorem}

\textsc{Remark. }The conclusions of Theorem \ref{T : Comparison} are
less precise than those of Theorems \ref{T : cumHRasy}--\ref{P :
path variance}, in the sense that they only apply to subsequences
$T_{n^{\prime }}$. Of course, this is due to the fact that, in the
statement of Theorem \ref{T : Comparison}, we do not make
\textit{any }assumption on the analytic properties of $k$ and $\nu
$, besides the conditions $\overline{k}\leq k\leq
\overline{\overline{k}}$ and $\overline{\nu }\leq \nu \leq \overline{%
\overline{\nu }}$. As will become clear in the subsequent sections,
more exact information can be deduced by adding some specific
requirements to the structure of $k$ and $\nu $.

\section{Applications\label{S : Applications}}

We will now consider noteworthy examples of random hazard rates by
specifying suitable kernels and the form of the background driving
CRM. In the following we will always consider CRMs with $\lambda $
being the Lebesgue measure on $\mathds{R}^{+}$, which appears a
natural choice in our context. This implies that Assumption (H1) is
met. Paragraph~\ref{linear} is
devoted to the study of the asymptotic behavior of the cumulative hazard $%
\tilde{H}$, whereas in Paragraph~\ref{quadratic} we deal with
quadratical functionals of the hazard rate.

\subsection{Asymptotics for the cumulative hazard\label{linear}}

As an illustration of Theorem \ref{T : cumHRasy}, we consider
different kernels and show how they are responsible for the rate of
divergence of the cumulative hazard and how they influence the
variance of the limiting Gaussian random variable in the CLT
(\ref{aaa}). We first consider general
homogeneous CRM such that $\int_{[1,\infty )}v^{4}\rho (\mathrm{d}v)<\infty $%
, which is tantamount of requiring the part of condition (H2)
involving the jump component of the Poisson intensity to be
satisfied. Moreover, set, for notational convenience, $K_{\rho
}^{\left( i\right) }=\int_{0}^{\infty
}s^{2}\rho (\mathrm{d}s)$, $i=1,2$, and $I_{i}=I_{i}\left( T\right) =\int_{%
\mathds{R}^{+}\times \mathbb{X}}\left[ k_{T}^{\left( 0\right)
}\left(
s,x\right) \right] ^{i}\nu \left( \mathrm{d}s,\mathrm{d}x\right) $ for $%
i=1,2,3$. Note that $I_{1}\left( T\right)
=\mathds{E}[\tilde{H}(T)]$.

\smallskip

\noindent \textsl{(i) Rectangular kernel.} The kernel $k\left( t,x\right) =%
\mathbb{I}_{\left( |t-x|\leq \tau \right) }$ where $\tau >0$
represents a bandwidth, is known as uniform rectangular kernel. Such
a kernel represents a sensible choice when no prior information on
the shape of the hazard rate is available. See, e.g., Ishwaran and
James (2004). In this setup (H2) is clearly met,
\[
k_{T}^{(0)}(s,x)=\left\{
\begin{array}{ll}
s\,(x+\tau ) & \qquad 0<x<\tau \\
s\,2\tau & \qquad \tau \leq x<T-\tau \\
s\,[T+\tau -x] & \qquad T-\tau \leq x<T+\tau \\
0 & \qquad \hbox{elsewhere}%
\end{array}%
\right.
\]%
and $k_{T}^{(0)}(s,x)\in L^{3}(\nu )$ for all $T>0$. We also have, as $%
T\rightarrow +\infty $, $I_{1}\left( T\right) $ $=K_{\rho }^{\left(
1\right) }\left\{ 2T\tau -\frac{1}{2}\tau ^{2}\right\} $ $=2\tau
K_{\rho }^{\left( 1\right) }T+o\left( T^{1/2}\right) $, $I_{2}\left(
T\right) \sim 4K_{\rho
}^{\left( 2\right) }\tau ^{2}T$ and $I_{3}\left( T\right) \sim cT$ for some $%
c>0$. Hence, \eqref{a-} and \eqref{aa} are satisfied with $%
C_{0}(k,T)=T^{-1/2}$ and, by Theorem \ref{T : cumHRasy}, we obtain
\begin{equation}
\frac{1}{\sqrt{T}}\left[ \tilde{H}(T)-2\tau K_{\rho }^{\left( 1\right) }T%
\right] \overset{\mathrm{law}}{\longrightarrow}X,
\label{eq:mean_rectangular}
\end{equation}%
where $X\sim \mathscr{N}\left( 0,4K_{\rho }^{\left( 2\right) }\tau
^{2}\right) $.

\smallskip

\noindent \textsl{(ii) Dykstra--Laud kernel.} If $k\left( t,x\right) =%
\mathbb{I}_{\left( 0\leq x\leq t\right) }$, then the random hazard
rate is monotone increasing. Such a kernel, which is widely
exploited in practice, was first proposed in Dykstra and Laud
(1981). It is easy to see that (H2) is satisfied and that
$k_{T}^{(0)}(s,x)=s(T-x)\mathbb{I}_{\left( 0\leq x\leq
T\right) }\in L^{3}(\nu )$ for all $T>0$. Moreover, one obtains $I_{1}=\frac{%
K_{\rho }^{\left( 1\right) }}{2}T^{2}$, $I_{2}=\frac{K_{\rho }^{(2)}}{3}%
T^{3} $ and $I_{3}=\frac{K_{\rho }^{(2)}}{4}T^{4}$, so that \eqref{a-} and %
\eqref{aa} are met with $C_{0}(k,T)=T^{-3/2}$. Hence, by Theorem
\ref{T : cumHRasy}, we have
\begin{equation}
\frac{1}{T^{\frac{3}{2}}}\left[ \tilde{H}(T)-\frac{K_{\rho }^{\left(
1\right) }}{2}T^{2}\right] \overset{\mathrm{law}}{\longrightarrow}X,
\label{eq:mean_DL}
\end{equation}%
where $X\sim \mathscr{N}\left( 0,\frac{K_{\rho }^{(2)}}{3}\right) $.
Note that the Dykstra-Laud cumulative hazard has a quadratic
asymptotic trend, whereas the trend obtained from a rectangular
kernel is linear. Moreover, the speed at which the Dykstra-Laud
cumulative hazard diverges from its trend is significantly faster
than in the rectangular case. The reason may be that the former
produces monotone increasing hazard rates whereas the latter not.
This phenomenon, well exemplified by our result, should be taken
into account when deciding which kernel to adopt.

\smallskip

\noindent \textsl{(iii) Ornstein--Uhlenbeck kernel.} If $k\left( t,x\right) =%
\sqrt{2\kappa }\exp \left( -\kappa \left( t-x\right) \right) \mathbb{I}%
_{\left( 0\leq x\leq t\right) }$, then the random hazard rate is an
Ornstein--Uhlenbeck--type process. Such models for the hazard rate
are employed in Nieto--Barajas and Walker (2004, 2005). In this
case, (H2) is
met, $k_{T}^{(0)}(s,x)=s\sqrt{2/\kappa }\,(1-\mathrm{e}^{-\kappa (T-x)})%
\mathbb{I}_{\left( 0\leq x\leq T\right) }\in L^{3}(\nu )$ for all
$T>0$, and we have that, as $T$ diverges to infinity, $I_{1}\left(
T\right) =$ $K_{\rho }^{\left( 1\right) }\sqrt{2/\kappa }\left\{
T-e^{-T}/\kappa +\kappa ^{-1}\right\} =K_{\rho }^{\left( 1\right)
}\sqrt{2/\kappa }T+o\left(
T^{1/2}\right) $, $I_{2}\left( T\right) \sim \frac{2K_{\rho }^{(2)}}{\kappa }%
T$ and $I_{3}\left( T\right) \sim cT$ for some constant $c>0$. Hence, %
\eqref{a-} and \eqref{aa} are satisfied with
$C_{0}(k,T)=T^{-1/2}$. From Theorem \ref{T : cumHRasy} it follows
that
\begin{equation}
\frac{1}{\sqrt{T}}\left[ \tilde{H}(T)-K_{\rho}^{(1)}\sqrt{\frac{2}{\kappa }}T%
\right] \overset{\mathrm{law}}{\longrightarrow}X,
\label{eq:mean_OU}
\end{equation}%
where $X\sim \mathscr{N}\left( 0,\frac{2K_{\rho }^{(2)}}{\kappa
}\right) $. One may note that the trend and the rate of divergence
from the trend associated with the Ornstein--Uhlenbeck kernel
coincide with those arising from the rectangular kernel. Moreover,
given the same background driving
CRM, the variances of the limiting Gaussian random variables appearing in (%
\ref{eq:mean_rectangular}) and (\ref{eq:mean_OU}) coincide if the
parameters are chosen in such a way that $\kappa =1/(2\tau ^{2})$.

\smallskip

\noindent \textsl{(iv) U-shaped or bath--tube kernel.} If $k\left(
t,x\right) =\mathbb{I}_{\left( |t-\beta |\geq x\right) }$ with
$\beta >0$, then the corresponding hazard rates are U--shaped with
minimum at $\beta $. Such a kernel is suggested by Lo and Wong
(1989). See also James (2003) and Ishwaran and James (2004). It is
easy to check that (H2) is met,
\[
k_{T}^{(0)}(s,x)=\left\{
\begin{array}{ll}
s\,(T-2x) & \qquad 0<x<\beta \\
s\,[T-(\beta +x)] & \qquad \beta \leq x<T-\beta \\
0 & \qquad \hbox{elsewhere}%
\end{array}%
\right.
\]%
and $k_{T}^{(0)}(s,x)\in L^{3}(\nu )$ for all $T>0$. Moreover, as $%
T\rightarrow +\infty $, $I_{1}\left( T\right) =\frac{1}{2}K_{\rho
}^{\left( 1\right) }T^{2}+o\left( T^{3/2}\right) $, $I_{2}\sim
\frac{K_{\rho }^{\left( 2\right) }}{3}T^{3}$ and $I_{3}\sim cT^{4}$
for some constant $c>0$. Choosing $C_{0}(k,T)=T^{-3/2}$, \eqref{a-}
and \eqref{aa} are satisfied and from Theorem \ref{T : cumHRasy} we
deduce
\begin{equation}
\frac{1}{T^{\frac{3}{2}}}\left[ \tilde{H}(T)-\frac{1}{2}K_{\rho
}^{\left( 1\right) }T^{2}\right]
\overset{\mathrm{law}}{\longrightarrow}X, \label{eq:mean_U_shaped}
\end{equation}%
where $X\sim \mathscr{N}\left( 0,\frac{K_{\rho }^{\left( 2\right) }}{3}%
\right) $. Note that the bath--tube kernel produces the same
asymptotic behaviour of the Dykstra and Laud kernel: this fact is
not surprising since after reaching its minimum in $\beta $, also
the bath--tube kernel is monotone increasing. Of course, one can
regard the Dykstra and Laud kernel as a degenerate bath--tube
kernel, corresponding to the case $\beta =0$.

As apparent from the statement of Theorem \ref{T : cumHRasy} and
from the discussion provided above, the variances of the limiting
Gaussian random
variables appearing in (\ref{aaa}), (\ref{eq:mean_rectangular}), (\ref%
{eq:mean_DL}), (\ref{eq:mean_OU}) and (\ref{eq:mean_U_shaped}),
always
depend on the jump part of the Poisson intensity. For instance, if $\tilde{%
\mu}$ is the generalized gamma CRM with intensity \eqref{eq:gg}, then $%
K_{\rho }^{(2)}=\frac{(1-\sigma )}{\gamma ^{2-\sigma }}$. This
confirms the empirical finding, used in tuning the prior parameters,
that a small $\gamma $ induces a large variance. To avoid confusion,
note that in the setting of e.g. Ishwaran and James (2004) $\beta
=1/\gamma $ and, hence, their claim that a large $\beta $ induces a
non--informative prior is coherent with our result. As for $\sigma
$, the variance is maximal in $\sigma =0$ if $\gamma \leq
\mathrm{e}$, whereas it is maximized in $\sigma =(\log (\gamma
)-1)/\log (\gamma )$ if $\gamma \geq \mathrm{e}$.

Let us now turn attention to hazards based on non--homogeneous CRM,
specifically the extended gamma and beta CRMs presented in Section 2. From %
\eqref{eq:extended_gamma} and \eqref{eq:beta} one can see that their
non--homogeneity is due to the strictly positive functions $\beta $
and $c$, respectively. According to their structure we distinguish
three cases: (a) if $\beta (x)=\bar{\beta}$ in
\eqref{eq:extended_gamma} and $c(x)=\bar{c}$ in \eqref{eq:beta}, the
CRMs become homogeneous and the previous results
hold with $K_{\rho }^{(2)}$ equal to $1/\bar{\beta}^{2}$ and $1/(1+\bar{c})$%
, respectively. (b) If $\beta $ (or $c$) are bounded by some finite
constant $M$, then one can apply Theorem 4 to conclude that
$C_{0}(k,T)$ has the same order as in the examples above, thus
depending on the choice of the kernel. Moreover, if $\beta $ (or
$c$) are eventually non--decreasing (non--increasing) the
convergence holds for any diverging sequence $T_{n}$ with the
variance of the limiting Gaussian random variable depending on the
choice of $\beta $ (or $c$) taking value in the range $[\sigma _{0}^{2}(%
\overline{k}),\sigma _{0}^{2}(\overline{\overline{k}})]$. (c) If
$\beta $ (or $c$) diverge to $+\infty $ as $x\rightarrow +\infty $,
quite interesting phenomena appear, which shed some light on the
possible use of the factor of
non--homogeneity represented by the functions $\beta $ (or $c$). Set, for $%
i=1,2,3$, $K_{\rho }^{(i)}(x)=\int_{0}^{\infty }s^{i}\rho
(\mathrm{d}s|x)$, so that $I_{i}$ becomes $\int_{\mathbb{X}}K_{\rho
}^{(i)}(x)\left[ \int_{0}^{T}k(t,x)\mathrm{d}t\right]
^{i}\mathrm{d}x$. For both CRMs, a diverging $\beta $ (or $c$)
implies that $K_{\rho }^{(2)}(x)\rightarrow 0$:
this, indeed, affects the asymptotic behavior of the cumulative hazard $%
\tilde{H}$. To be more specific, consider the Dykstra and Laud
kernel
combined with an extended gamma CRM such that $\beta (x)\sim \sqrt{x}$ as $%
x\rightarrow \infty $: it follows that $I_{2}\sim \log (T)T^{2}$ and $%
I_{3}\sim dT^{3}$ for some constant $d>0$. Hence, \eqref{a-} and
\eqref{aa}
are satisfied with $C_{0}(k,T)=(\sqrt{\log (T)}\,T)^{-1}$ and, by Theorem %
\ref{T : cumHRasy}, we have
\begin{equation}
\frac{1}{\sqrt{\log (T)}\,T}\left[ \tilde{H}(T)-\mathds{E}[\tilde{H}(T)]%
\right] \overset{\mathrm{law}}{\longrightarrow }X,
\label{eq:mean_extgamma1}
\end{equation}%
where $X\sim \mathscr{N}\left( 0,1\right) $. Comparing %
\eqref{eq:mean_extgamma1} with \eqref{eq:mean_DL} one notes that the
rate of divergence from the trend $\mathds{E}[\tilde{H}(T)]$ is
reduced from $T^{3/2} $ to $\sqrt{\log (T)}\,T$. As for
$\mathds{E}[\tilde{H}(T)]$, it is important to remark that the
overall growth (though not the dominating term which is
$4/3\,T^{3/2}$) depends on the particular form of $\beta $. Still
assuming $\beta (x)\sim \sqrt{x}$ and letting $b$ be a positive
constant, we
obtain, for instance, $\mathds{E}[\tilde{H}(T)]=4/3T^{3/2}+o(T\sqrt{\log (T)}%
)$ when $\beta (x)=\mathbb{I}_{(0,b]}(x)+x^{1/2}\mathbb{I}_{(b,\infty )}(x)$%
, and $\mathds{E}[\tilde{H}(T)]=4/3\,T^{3/2}-\log (T)T+o(T\sqrt{\log
(T)})$
if $\beta (x)=(1+x^{1/2})$. Again, comparing these findings with %
\eqref{eq:mean_DL} it is apparent that the trend has been reduced
from $T^{2}
$ to $T^{3/2}+o(T^{3/2})$. On the other hand, with the beta CRM, we have $%
K_{\rho }^{(1)}(x)=1$ and, consequently, $I_{1}(T)=\mathds{E}[\tilde{H}(T)]%
=1/2\,T^{2}$ whatever the choice of $c$. Selecting $c(x)\sim \sqrt{x}$ as $%
x\rightarrow \infty $, we obtain $I_{2}\sim 16/15\,T^{5/2}$ and
$I_{3}\sim
d\log (T)T^{3}$ for some constant $d>0$. Thus, with $C_{0}(k,T)=T^{-5/4}$, %
\eqref{a-} and \eqref{aa} are met and Theorem \ref{T : cumHRasy} yields%
\begin{equation}
\frac{1}{T^{\frac{5}{4}}}\left[
\tilde{H}(T)-\frac{1}{2}\,T^{2}\right]
\overset{\mathrm{law}}{\longrightarrow }X,  \nonumber
\end{equation}%
where $X\sim \mathscr{N}\left( 0,16/15\right) $. Hence, compared
with the homogeneous case in \eqref{eq:mean_DL}, the beta CRM does
not affect the trend but still decreases the rate of divergence from
$T^{2}$ to $T^{5/4}$.

\noindent If, instead, we consider the rectangular kernel with
$\tau=1$ combined with an extended gamma CRM such that again
$\beta(x)\sim \sqrt{x}$ as $x\rightarrow \infty$, it follows that
$I_2\sim 4 \log(T)$ and $I_3\to d$ for some constant $d>0$. Hence,
\eqref{a-} and \eqref{aa} are satisfied with
$C_0(k,T)=(\sqrt{\log(T)})^{-1}$ and, by Theorem \ref{T : cumHRasy},
we have
\begin{equation}
\frac{1}{\sqrt{\log(T)}}\left[ \tilde{H}(T)-\mathds{E}[\tilde{H}(T)]
\right] \overset{\mathrm{law}}{\longrightarrow}X,  \nonumber
\end{equation}
where $X\sim \mathscr{N}\left( 0, 4\right)$. Hence, we see that the
rate of divergence from $\mathds{E}[\tilde{H}(T)]$ has been reduced
with respect to the homogeneous case in \eqref{eq:mean_rectangular}
decreasing from $T^{1/2}$ to $\sqrt{\log(T)}$. As before,
$I_1(T)=\mathds{E}[\tilde{H}(T)]$ depends on the particular form of
$\beta$. With $\beta (x)\sim \sqrt{x}$, $b$ being a
positive constant, we have $I_1(T)=4 T^{1/2}+o(\sqrt{\log(T)})$ if $\beta(x)=%
\mathbb{I}_{(0,b]}(x)+x^{1/2}\mathbb{I}_{(b,\infty)}(x)$ and
$I_1(T)=4 T^{1/2}-2\log(T)+o(\sqrt{\log(T)})$ if
$\beta(x)=(1+x^{1/2})$. By comparing these trends with the one in
\eqref{eq:mean_rectangular} one can appreciate its reduction from
$T$ to $T^{1/2}+o(T^{1/2})$.

Replacing the extended gamma CRM with a beta process we have
$I_1(T)=2\,
T-1/2$ whatever the choice of $c$. Moreover, if $c(x)\sim \sqrt{x}$ as $%
x\rightarrow \infty $ we obtain $I_2\sim 8 \, T^{1/2}$ and $I_3\sim
d
\log(T) $ for some $d>0 $. By setting $C_0(k,T)=T^{-1/4}$ \eqref{a-} and %
\eqref{aa} are met and Theorem \ref{T : cumHRasy} leads to
\begin{equation}
\frac{1}{T^{\frac{1}{4}}}\left[ \tilde{H}(T)-2 \, T \right] \overset{\mathrm{%
law}}{\longrightarrow}X,  \nonumber
\end{equation}
where $X\sim \mathscr{N}\left( 0, 8\right)$. Hence, with respect to %
\eqref{eq:mean_rectangular}, the trend is unchanged and the rate of
divergence halved.

By means of the previous examples the impact of a non--homogeneous
CRM becomes apparent: a non--homogeneous CRM allows to reduce both
the trend of the cumulative hazard and the rate at which it diverges
from its trend. An extended gamma CRM is able to reduce both,
whereas a beta CRM affects only the rate of divergence from the
trend. Overall, by studying also other examples, not reported here,
of functions $\beta $ and $c $ with the 4 different kernels
considered above, some interesting indications can be drawn. For
instance, denote by $T^{\eta }$ the rate at which the cumulative
hazard based on the homogeneous version of an extended gamma (or
beta) CRM diverges from its trend (e.g. $\eta =3/2$ in the
Dykstra-Laud case). Then, by choosing a suitable diverging $\beta $
(or $c$) the rate can be tuned at any order in the range $[T^{\eta
-1/2},T^{\eta }]$. Analogous conclusions can be derived for the
trend when using a hazard based on an extended gamma
CRM: the trend corresponding to the homogeneous case $T^{\alpha}$ (e.g. $%
\alpha=2$ for the Dykstra--Laud kernel) can be tuned by the choice
of $\beta$ at any rate in the range $[T^{\alpha-1}, T^{\alpha}]$.

\subsection{Asymptotics for quadratic functionals\label{quadratic}}

In this paragraph we consider quadratic functionals of the random
hazard rate. We derive central limit theorems for the path--second
moments and the path--variances of hazard rates with specific
kernels and driving CRM. Our results will be mainly based on
Theorems~\ref{T : path2ndM} and \ref{P : path variance}. As in the
previous paragraph, we first deal with general
homogeneous CRM such that $\int_{[1,\infty )}v^{4}\rho (\mathrm{d}v)<\infty $%
; this requirement combined with the structure of kernels we
consider ensures that (H2) is satisfied. Finally set, as before,
$K_{\rho }^{(i)}=\int_{0}^{\infty }s^{i}\rho (\mathrm{d}s)$, for
$i=1,2,3,4$.
\smallskip

\noindent \textsl{(i) Rectangular kernel.} We start by considering
the rectangular kernel and derive CLTs for the path--second moment
and for the
path--variance of hazard rates. Some simple calculations lead to write, for $%
x>y$ and $T>2\tau $,
\[
k_{T}^{(1)}(s,x;t,y)=\left\{
\begin{array}{ll}
\frac{st}{T}\,(y+\tau ) & y<x<\tau ,0<y<\tau  \\
\frac{st}{T}\,(y+2\tau -x) & (\tau \vee y)\leq x<(y+2\tau ),0<y<T-\tau  \\
\frac{st}{T}\,[T+\tau -x] & y\leq x<T+\tau ,\ T-\tau \leq y<T+\tau  \\
0 & \hbox{elsewhere}%
\end{array}%
\right.
\]%
Moreover, $k_{T}^{(2)}(s,x)=s\,T^{-1}k_{T}^{(0)}(s,x)$ and for
$T>2\tau $, one has
\[
k_{T}^{(3)}(s,x)=\left\{
\begin{array}{ll}
\frac{sK_{\rho }^{(1)}}{T}\,\left[ \frac{1}{2}x^{2}+\tau x\right]  &
0<x<\tau  \\
\frac{sK_{\rho }^{(1)}}{T}\,\left[ -\frac{1}{2}x^{2}+2\tau x\right]
& \tau
\leq x<2\tau  \\
\frac{sK_{\rho }^{(1)}}{T}\,2\tau ^{2} & 2\tau \leq x<T-\tau  \\
\frac{sK_{\rho }^{(1)}}{T}\,\left[ \!-\frac{1}{2}T^{2}\!+\!T(x+\tau )\!+\!%
\frac{3}{2}\tau ^{2}\!\!-\!\tau x\!-\!\frac{1}{2}x^{2}\!\right]  &
T-\tau
\leq x<T+\tau  \\
0 & \hbox{elsewhere}%
\end{array}%
\right.
\]%
In order to apply Theorem~\ref{T : path2ndM} let us first consider
Condition
1., which allows to determine the rate function: it turns out that $%
C_{1}\left( k,T\right) =\sqrt{T}$ since
\begin{equation}
2\,T\,\left\Vert k_{T}^{\left( 1\right) }\right\Vert _{L^{2}\left(
\nu ^{2}\right) }^{2}\rightarrow \sigma _{1}^{2}(k)=\frac{16\,\tau
^{3}\,(K_{\rho }^{(2)})^{2}}{3}  \label{-w}
\end{equation}%
The verification of Conditions 2.--6. can be achieved by simple
though quite lengthy calculations.

Indeed, letting, for $i=1,\ldots ,4$, $d_{i}$ be a positive
constant, one obtains
\begin{align*}
& \hbox{2.}\ \ T^{2}\,\left\Vert k_{T}^{\left( 1\right) }\right\Vert
_{L^{4}\left( \nu ^{2}\right) }^{4}\sim \frac{d_{1}}{T}\rightarrow 0 \\
& \hbox{3.}\ \ T^{2}\,\left\Vert k_{T}^{\left( 1\right) }\star
_{1}^{1}k_{T}^{\left( 1\right) }\right\Vert _{L^{2}\left( \nu
^{2}\right)
}^{2}\sim \frac{d_{2}}{T}\rightarrow 0 \\
& \hbox{4.}\ \ T^{2}\,\left\Vert k_{T}^{\left( 1\right) }\star
_{2}^{1}k_{T}^{\left( 1\right) }\right\Vert _{L^{2}\left( \nu
\right)
}^{2}\sim \frac{d_{3}}{T}\rightarrow 0 \\
& \hbox{5.}\ \ T\left\Vert k_{T}^{\left( 2\right)
}\!\!+\!2k_{T}^{\left( 3\right) }\right\Vert _{L^{2}\left( \nu
\right) }^{2}\!\!\!\!\!\rightarrow \sigma _{2}^{2}(k)=16\tau
^{2}\!\left[ \frac{K_{\rho }^{(4)}}{4}+\tau K_{\rho }^{(3)}K_{\rho
}^{(1)}+\tau ^{2}K_{\rho }^{(2)}\left( K_{\rho
}^{(1)}\right) ^{2}\right] \\
& \hbox{6.}\ \ T^{^{\frac{3}{2}}}\,\left\Vert k_{T}^{\left( 2\right)
}+2k_{T}^{\left( 3\right) }\right\Vert _{L^{3}\left( \nu \right)
}^{3}\sim \frac{d_{4}}{T^{1/2}}\rightarrow 0
\end{align*}%
Since
\begin{equation}
\frac{1}{T}\int_{0}^{T}\mathds{E}(\tilde{h}(t)^{2})\mathrm{d}t=2\tau
\,K_{\rho }^{\left( 2\right) }+4\,\tau ^{2}\,\left( K_{\rho
}^{\left( 1\right) }\right) ^{2}+o\left( T^{-1/2}\right) ,
\label{eq:h2asymp}
\end{equation}%
we deduce from Theorem \ref{T : path2ndM} the following asymptotic
result, concerning the path-second moment of $\tilde{h}(t)$:
\begin{equation}
T^{1/2}\left\{
\frac{1}{T}\int_{0}^{T}\tilde{h}(t)^{2}\mathrm{d}t-\left( 2\tau
\,K_{\rho }^{\left( 2\right) }+4\,\tau ^{2}\,\left( K_{\rho
}^{\left(
1\right) }\right) ^{2}\right) \right\} \overset{\mathrm{law}}{%
\longrightarrow }X\text{,}  \nonumber
\end{equation}%
where $X\sim \mathscr{N}\left( 0,\sigma _{1}^{2}(k)+\sigma
_{2}^{2}(k)\right) $ with
\[
\sigma _{1}^{2}(k)+\sigma _{2}^{2}(k)=16\tau ^{2}\left[ \frac{K_{\rho }^{(4)}%
}{4}+\tau K_{\rho }^{(3)}K_{\rho }^{(1)}+\frac{\tau \left( K_{\rho
}^{(2)}\right) ^{2}}{3}+\tau ^{2}K_{\rho }^{(2)}\left( K_{\rho
}^{(1)}\right) ^{2}\right] .
\]%
Now we concentrate on a CLT involving the path--variance of
$\tilde{h}(t)$, that we shall obtain as an application of Theorem
\ref{P : path variance}. In particular, we must verify that
Conditions 1, 2 and 3 in the statement of such result are verified,
for some appropriate positive constants $\delta \left( k\right) $
and $\sigma _{3}^{2}\left( k\right) $. Indeed, one has that, as
$T\rightarrow +\infty $,
\begin{align}
\frac{C_{1}\left( k,T\right) }{\left( TC_{0}\left( k,T\right)
\right) ^{2}}&
=T^{-\frac{1}{2}}\rightarrow 0  \label{w} \\
\frac{2C_{1}\left( k,T\right) }{T^{2}C_{0}\left( k,T\right) }\, \mathds{E}%
\lbrack \tilde{H}(T)]& =\frac{2}{T}\left\{ 2\tau K_{\rho }^{\left(
1\right) }T+o\left( T\right) \right\} \rightarrow 4\tau K_{\rho
}^{\left( 1\right) }:=\delta \left( k\right)  \label{ww}
\end{align}%
and also
\begin{multline}
\left\Vert \ C_{1}\left( k,T\right) \left( k_{T}^{\left( 2\right)
}+2k_{T}^{\left( 3\right) }\right) -\delta \left( k\right)
C_{0}\left( k,T\right) k_{T}^{\left( 0\right) }\right\Vert
_{L^{2}\left( \nu \right)
}^{2} \\
\rightarrow 16\tau ^{2}\left[ \frac{K_{\rho }^{(4)}}{4}-\tau K_{\rho
}^{(3)}K_{\rho }^{(1)}+\tau ^{2}K_{\rho }^{(2)}\left( K_{\rho
}^{(1)}\right) ^{2}\right] :=\sigma _{3}^{2}\left( k\right) .
\label{www}
\end{multline}%
\noindent The fact that $\mathds{E}\lbrack \tilde{H}(T)]=K_{\rho
}^{\left(
1\right) }\left\{ 2T\tau -\frac{1}{2}\tau ^{2}\right\} $ combined with %
\eqref{eq:h2asymp} yields
\[
\frac{1}{T}\int_{0}^{T}\mathds{E}\left[
\tilde{h}(t)-\frac{\mathds{E}\lbrack \tilde{H}(T)]}{T}\right]
^{2}\mathrm{d}t=2\tau K_{\rho }^{(2)}+o(T^{-1/2}).
\]%
Hence, by using (\ref{w})--(\ref{www}), we deduce from Theorem
\ref{P : path variance} that
\begin{equation}
\sqrt{T}\times \left\{ \frac{1}{T}\int_{0}^{T}[\tilde{h}(t)-\frac{1}{T}%
\tilde{H}(T)]^{2}\mathrm{d}t-2\tau K_{\rho }^{(2)}\right\} \overset{\mathrm{%
law}}{\longrightarrow }X,  \nonumber
\end{equation}%
where $X\sim \mathscr{N}\left( 0,\sigma _{1}^{2}\left( k\right)
+\sigma _{3}^{2}\left( k\right) \right) $, and $\sigma
_{1}^{2}\left( k\right) $ and $\sigma _{3}^{2}\left( k\right) $ are
given by (\ref{-w}) and (\ref{www}), respectively.

\smallskip

\noindent \textsl{(ii) Ornstein--Uhlenbeck kernel.} Let us now
derive the CLT for the path--second moment and the path--variance of
hazards based on the Ornstein--Uhlenbeck kernel. For this case we
easily obtain
\begin{align*}
k_{T}^{(1)}(s,x;t,y)& =\frac{s\,t}{T}\mathrm{e}^{\kappa (x+y)}\left( \mathrm{%
e}^{-2\kappa x}-\mathrm{e}^{-2\kappa T}\right) \mathbb{I}_{(0\leq
y\leq
x\leq T)} \\
k_{T}^{(2)}(s,x)& =\frac{s^{2}}{T}\mathrm{e}^{2\kappa x}\left( \mathrm{e}%
^{-2\kappa x}-\mathrm{e}^{-2\kappa T}\right) \mathbb{I}_{(0\leq x\leq T)} \\
k_{T}^{(3)}(s,x)& =\frac{sK_{\rho }^{(1)}}{\kappa \,T}\left( \mathrm{e}%
^{-2\kappa T}-\mathrm{e}^{-2\kappa x}\right) \,\left( \mathrm{e}^{\kappa x}-%
\mathrm{e}^{2\kappa x}\right) \mathbb{I}_{(0\leq x\leq T)}
\end{align*}%
and some tedious algebra allows to derive also $k_{T}^{(1)}\star
_{1}^{1}k_{T}^{(1)}$ and $k_{T}^{(1)}\star _{2}^{1}k_{T}^{(1)}$.
Condition 1. in Theorem~\ref{T : path2ndM} is verified by choosing
$C_{1}\left( k,T\right) =\sqrt{T}$: indeed,
\begin{equation}
2\,T\,\left\Vert k_{T}^{\left( 1\right) }\right\Vert _{L^{2}\left(
\nu
^{2}\right) }^{2}\rightarrow \sigma _{1}^{2}(k)=\frac{(K_{\rho }^{(2)})^{2}}{%
\kappa }.  \label{eq:s1OU}
\end{equation}%
Standard calculations allow to verify the validity of the other
conditions
in the statement of Theorem~\ref{T : path2ndM}. In particular, by letting $%
d_{i}$ ($i=1,\ldots ,4$) be a positive constant, one obtains
\begin{align*}
& \hbox{2.}\quad T^{2}\left\Vert k_{T}^{\left( 1\right) }\right\Vert
_{L^{4}\left( \nu ^{2}\right) }^{4}\sim \frac{d_{1}}{T}\rightarrow 0 \\
& \hbox{3.}\quad T^{2}\left\Vert k_{T}^{\left( 1\right) }\star
_{1}^{1}k_{T}^{\left( 1\right) }\right\Vert _{L^{2}\left( \nu
^{2}\right)
}^{2}\sim \frac{d_{2}}{T}\rightarrow 0 \\
& \hbox{4.}\quad T^{2}\left\Vert k_{T}^{\left( 1\right) }\star
_{2}^{1}k_{T}^{\left( 1\right) }\right\Vert _{L^{2}\left( \nu
\right)
}^{2}\sim \frac{d_{3}}{T}\rightarrow 0 \\
& \hbox{5.}\quad T\left\Vert k_{T}^{\left( 2\right) }+2k_{T}^{\left(
3\right) }\right\Vert _{L^{2}\left( \nu \right) }^{2}\rightarrow
\sigma
_{2}^{2}(k)=K_{\rho }^{(4)}+\frac{4}{\kappa }K_{\rho }^{(3)}K_{\rho }^{(1)}+%
\frac{4}{\kappa ^{2}}K_{\rho }^{(2)}\left( K_{\rho }^{(1)}\right) ^{2} \\
& \hbox{6.}\quad T^{^{\frac{3}{2}}}\left\Vert k_{T}^{\left( 2\right)
}+2k_{T}^{\left( 3\right) }\right\Vert _{L^{3}\left( \nu \right)
}^{3}\sim \frac{d_{4}}{T^{1/2}}\rightarrow 0
\end{align*}%
Since, as $T\rightarrow +\infty $,
\begin{equation}
\frac{1}{T}\int_{0}^{T}\mathds{E}(\tilde{h}(t)^{2})\mathrm{d}t=K_{\rho
}^{(2)}+\frac{2\,\left( K_{\rho }^{(1)}\right) ^{2}}{\kappa
}+o\left( T^{-1/2}\right) ,  \label{eq:OUh2}
\end{equation}%
we deduce from Theorem~\ref{T : path2ndM} the following result for
the path--second moment:
\begin{equation}
T^{1/2}\left\{
\frac{1}{T}\int_{0}^{T}\tilde{h}(t)^{2}\mathrm{d}t-\left[ K_{\rho
}^{(2)}+\frac{2\,\left( K_{\rho }^{(1)}\right) ^{2}}{\kappa }\right]
\right\} \overset{\mathrm{law}}{\longrightarrow }X\text{,}
\nonumber
\end{equation}%
where $X\sim \mathscr{N}\left( 0,K_{\rho }^{(4)}+\frac{4}{\kappa
}K_{\rho
}^{(3)}K_{\rho }^{(1)}+\frac{(K_{\rho }^{(2)})^{2}}{\kappa }+\frac{4}{k^{2}}%
K_{\rho }^{(2)}\left( K_{\rho }^{(1)}\right) ^{2}\right) $. As far
as the path--variance is concerned, one verifies easily that the
conditions of
Theorem \ref{P : path variance} are verified, with $\delta \left( k\right) =%
\frac{2^{3/2}}{\sqrt{\kappa }}K_{\rho }^{\left( 1\right) }$ and
\begin{equation}
\sigma _{3}^{2}\left( k\right) :=K_{\rho }^{\left( 4\right) }-\frac{4}{%
\kappa }\,K_{\rho }^{\left( 3\right) }\,K_{\rho }^{\left( 1\right) }+\frac{4%
}{\kappa ^{2}}\,\,K_{\rho }^{\left( 2\right) }\left( K_{\rho
}^{\left( 1\right) }\right) ^{2}  \label{u}
\end{equation}

\noindent Using \eqref{eq:OUh2}, it is straightforward to see that
\[
\frac{1}{T}\int_{0}^{T}\mathds{E}\left[
\tilde{h}(t)-\frac{\mathds{E}\lbrack \tilde{H}(T)]}{T}\right]
^{2}\mathrm{d}t=K_{\rho }^{(2)}+o(T^{-1/2}).
\]%
As a consequence, we deduce from Theorem \ref{P : path variance}
that
\begin{equation}
\sqrt{T}\times \left\{ \frac{1}{T}\int_{0}^{T}[\tilde{h}(t)-\frac{1}{T}%
\tilde{H}(T)]^{2}\mathrm{d}t-K_{\rho }^{\left( 2\right) }\right\} \overset{%
\mathrm{law}}{\longrightarrow }X,  \nonumber
\end{equation}%
with $X\sim \mathscr{N}\left( 0,\sigma _{1}^{2}\left( k\right)
+\sigma _{3}^{2}\left( k\right) \right) $, where $\sigma
_{1}^{2}\left( k\right) $
and $\sigma _{3}^{2}\left( k\right) $ are given by \eqref{eq:s1OU} and %
\eqref{u}, respectively.

Before considering the Dykstra and Laud kernel and the U--shaped
kernel, let us make the previous results completely explicit by
specifying the background driving CRM. For both the rectangular and
the Ornstein--Uhlenbeck kernel the rate function is the same and the
CRM affects the variance of the limiting Gaussian random variable
for both path--second moment and path-variance of the hazard rate.
Take, as before the generalized gamma CRM
with Poisson intensity \eqref{eq:gg} and denote the Pochhammer symbol by $%
(a)_{n}:=\Gamma (a+n)/\Gamma (a)$. For this choice we have $K_{\rho
}^{(c)}=[(1-\sigma )_{c-1}](\gamma ^{c-\sigma })^{-1}$ for any
$c>0$. For the Ornstein--Uhlenbeck kernel the variance is then given
by
\begin{equation}
\sigma _{1}^{2}(k)+\sigma _{2}^{2}(k)=\frac{(1-\sigma )\,\left(
4\kappa ^{-1}\gamma ^{2\sigma }+(9-5\sigma )\gamma ^{\sigma }+\kappa
(2-\sigma )_{2}\right) }{\kappa \gamma ^{4-\sigma }}
\label{eq:gg2nd}
\end{equation}%
which decreases as $\kappa $ and $\gamma $ increase for any given
$(\gamma ,\sigma )$ and $(\kappa ,\sigma )$, respectively. Moreover,
it is maximized by $\sigma =0$ for low values of $\kappa $ and
$\gamma $, whereas, for moderately large values of $\kappa $ and
$\gamma $, the maximizing $\sigma $ increases as $\kappa $ and
$\gamma $ increase. For instance, if $\kappa =0.5$ and $\gamma =2$,
the maximizing $\sigma $ is approximately equal to $0.22$ and the
overall variance is $2.56$. To highlight the incidence of the prior
parameters note that with $\kappa =1$ and $\gamma =3$, the maximizing $%
\sigma $ and the variance are approximately equal to $0.52$ and
$0.29$, respectively. Using the asymptotic variance as a guideline
for fixing the prior parameters seems a sensible and straightforward
choice since it summarizes in a single expression the various
effects of the parameters. Turning to the path--variance a hazard
based on a generalized gamma CRM with Ornstein--Uhlenbeck kernel
will have variance given by
\begin{equation}
\sigma _{1}^{2}(k)+\sigma _{3}^{2}(k)=\frac{(1-\sigma )\,\left(
4\kappa ^{-1}\gamma ^{2\sigma }-(7-3\sigma )\gamma ^{\sigma }+\kappa
(2-\sigma )_{2}\right) }{\kappa \gamma ^{4-\sigma }}
\label{eq:ggvar}
\end{equation}%
which behaves in the same way as \eqref{eq:gg2nd} but, obviously,
leads to
smaller values. Considering the same set of parameters as above we have: if $%
\kappa =0.5$ and $\gamma =2$, $\sigma \approx 0.61$ maximizes %
\eqref{eq:ggvar} and its value is $0.92$; if $\kappa =1$ and $\gamma =3$, %
\eqref{eq:ggvar} is maximized by $\sigma \approx 0.76$ leading to a
variance of $0.09$. Similar considerations hold also for the
asymptotic variance of a hazard based on the rectangular kernel
combined with a generalized gamma CRM.

Turning attention to quadratic functionals of hazards based on
non--homogeneous CRM the importance of our Theorem~\ref{T :
Comparison} becomes apparent: the verification of the conditions of
Theorem~\ref{T : path2ndM} and \ref{P : path variance} become
extremely difficult if not impossible. Hence, when it is possible to
bound above and below the Poisson
intensity of a non--homogeneous CRM so to meet the conditions of Theorem \ref%
{T : Comparison}, we are still able to state that the rate function is $%
C_{1}(k,T)=T^{1/2}$ for hazards based on rectangular and
Ornstein--Uhlenbeck
kernels. Moreover, we can deduce the convergence, along some subsequence $%
T_{n^{\prime }}$ of every diverging sequence $T_{n}$, of the
path--second moment and of the path-variance to a Gaussian random
variable with variance
taking value in the range $[\sigma _{1}^{2}(\overline{k})+\sigma _{2}^{2}(%
\overline{k}),\sigma _{1}^{2}(\overline{\overline{k}})+\sigma _{2}^{2}(%
\overline{\overline{k}})]$ and $[\sigma
_{1}^{2}(\overline{k})+\sigma _{3}^{2}(\overline{k}),\sigma
_{1}^{2}(\overline{\overline{k}})+\sigma
_{3}^{2}(\overline{\overline{k}})]$, respectively. In order to
deduce convergence for every diverging sequence, the structure of
the Poisson intensity has to be specified as well. Thus, let us
consider again the
extended gamma and beta CRMs. As noted in Section 4.1, supposing $\beta (x)=%
\bar{\beta}$ in \eqref{eq:extended_gamma} and $c(x)=\bar{c}$ in %
\eqref{eq:beta}, the CRMs become homogeneous and the previous
results hold with the same rate functions. Note that, for $a>0$,
$K_{\rho }^{(a)}=\Gamma (a)\,\bar{\beta}^{-a}$ in the extended gamma
case and $K_{\rho }^{(a)}=\Gamma (a)\,[(1+\bar{c})_{a-1}]^{-1}$ in
the beta case. Hence, with an Ornstein--Uhlenbeck kernel the
asymptotic variance of the path--second moment is equal to
$(\bar{\beta}^{4}\,\kappa ^{2})^{-1}\,(6\kappa
^{2}+9\kappa +4)$ for the former and equal to $[\kappa ^{2}(1+\bar{c})\,(1+%
\bar{c})_{3}]^{-1}\,(9\kappa \bar{c}^{2}+37\kappa \bar{c}+30\kappa
+6\kappa ^{2}(1+\bar{c})+30\kappa (1+\bar{c})+4(1+\bar{c})_{3})$ for
the latter. For the path--variance similar expressions are obtained.
If $\beta $ (or $c$) are functions bounded by some finite constant
$M$, then we are in the genuinely non--homogeneous case and, as
mentioned above, by Theorem 4 CLTs along subsequences of diverging
sequences are granted. To achieve
convergence along any sequence, it is enough to suppose that $\beta $ (or $c$%
) are eventually non--decreasing (or non--increasing), which
represents a sensible choice in any application. For instance,
considering an extended gamma CRM with non--decreasing $\beta $
taking values in $[L,M]$ combined with an Ornstein--Uhlenbeck kernel
the path--second moment will converge, along any sequence, to a
Gaussian random variable with variance $\sigma _{1}^{2}(k)+\sigma
_{2}^{2}(k)=(M^{4}\,\kappa ^{2})^{-1}\,(6\kappa ^{2}+9\kappa +4)$.
Analogous considerations hold for the path--variance.

\smallskip

\noindent \textsl{(iii) Dykstra--Laud and U--shaped kernels.} Our
results for quadratic functionals do not apply when choosing the
kernel $k$ to be the Dykstra--Laud or U--shaped kernel. Indeed, for
both kernels Conditions 3., 5. and 6. in Theorem~\ref{T : path2ndM}
are not met. Moreover, also the additional conditions 1.--3. in
Theorem~\ref{P : path variance} are not satisfied. Note that
Condition 3. represents the most delicate since it involves a
contraction. Consider first the Dykstra--Laud kernel. It is easy to
see that $k_{T}^{(1)}(s,x;t,y)=\frac{s\,t}{T}(T-x)\mathbb{I}_{(0\leq
y\leq x\leq T)}$ and that $k_{T}^{(1)}\star
_{1}^{1}k_{T}^{(1)}\left(
s,x;t,y\right) =\frac{s\,tK_{\rho }\,(T-x)}{T^{2}}\left[ \frac{(T-y)^{2}}{2}-%
\frac{(T-x)^{2}}{2}\right] \mathbb{I}_{(0\leq y\leq x)}$. As for
Condition 1. we obtain with the choice $C_{1}=T^{-1}$
\[
\frac{2}{T^{2}}\left\Vert k_{T}^{\left( 1\right) }\right\Vert
_{L^{2}\left( \nu ^{2}\right) }^{2}\rightarrow \frac{K_{\rho
}^{2}}{6}.
\]%
This, however implies that the quantity in Condition 3. converges to
a
positive constant and the ones in Condition 5 and 6. diverge. In Theorem~\ref%
{P : path variance} we obtain that the quantity in Condition 1. is equal to $%
1$ and the one in Condition 2. diverges. Finally, Condition 3.
cannot be satisfied since Condition 5. in Theorem \ref{T : path2ndM}
is violated. For the U--shaped kernel we obtain again
$C_{1}(k,T)=T^{-1}$ and the asymptotic behaviour of the various
quantities involved in the conditions is the same as the one of the
Dykstra and Laud kernel. We have also tried with
non--homegeneous CRM: indeed, it seems possible to obtain $%
C_{1}(k,T)=T^{-\eta }$ with any $\eta \in (0,1]$, but the conditions
are nonetheless violated.

\noindent The fact that our results do not work for the
Dykstra--Laud and U--shaped kernels seem to suggest that kernels
yielding monotone increasing hazards (at least from some point
onwards as it is the case for the U--shaped kernel) exhibit a too
strong growth to be compatible with our conditions. Future research
will focus, on one side, on the translation of the conditions into
simple and intuitive sufficient ones regarding the behaviour of the
hazard rate induced by different classes of kernels and, on the
other side, to relax the conditions in order to cover models for
monotone increasing hazards.

\section{Proofs and further techniques\label{S : Proofs}}

In this section we collect the proofs of the main results of the
paper. As anticipated, we shall make a substantial use of the CLTs,
for sequences of single and double Poisson integrals, recently
established by Peccati and Taqqu (2006b). In the next subsection we
present some preliminary results concerning double Wiener-It\^{o}
integrals, with special attention devoted to weak convergence and
central limit theorems. Virtually all of the needed background
material, about stochastic integrals of any order with respect to
Poisson measures, can be found in Surgailis (1984) and in Chapter 10
of Kwapie\'{n} and Woyczy\'{n}ski (1992). A different approach,
based on Hilbert space techniques, is described in Nualart and Vives
(1990). The reader is also referred to Surgailis (2000) for an
updated review of related convergence results.

\subsection{Double integrals and CLTs\label{SS : 2ble int CLT}}

Throughout this section we consider a Poisson CRM $\tilde{N}$ such that (\ref%
{eq:H1}) is verified. Recall that $\tilde{N}^{c}$ is the compensated
Poisson measure defined in formulae (\ref{eq : Comp P}) and (\ref{eq
: Comp P2}).
For every $f\in L_{s,0}^{2}\left( \nu ^{2}\right) $, we denote by $I_{2}^{%
\tilde{N}^{c}}\left( f\right) $ the \textsl{double Wiener-It\^{o} integral }%
of $f$ with respect to $\tilde{N}^{c}$. The reader is referred to
Surgailis (1984) for precise definitions. Here, we shall recall
that, if $f\in L_{s,0}^{2}\left( \nu ^{2}\right) $ is a piecewise
constant function with support contained in a product set $S\times
S\subset (\mathds{R}^{+}\times \mathbb{X)}^{2}$ such that $\nu
\left( S\right) <+\infty $, then\ the variable
$I_{2}^{\tilde{N}^{c}}\left( f\right) $ is a genuine
(\textquotedblleft pathwise\textquotedblright ) double integral with
respect
to the restriction to $S\times S$ of the (signed) product measure $\tilde{N}%
^{c}\left( \mathrm{d}s,\mathrm{d}x\right) \tilde{N}^{c}\left( \mathrm{d}t,%
\mathrm{d}y\right) $. The very nature of $f$ implies that the
integration is
performed on the intersection between $S\times S$ and the non-diagonal set $%
D_{0}^{2}$. For a general $f\in L_{s,0}^{2}\left( \nu ^{2}\right) $, $I_{2}^{%
\tilde{N}^{c}}\left( f\right) $ is simply the limit in $L^{2}\left( %
\mathds{P}\right) $ of random variables of the kind $I_{2}^{\tilde{N}%
^{c}}\left( f_{k}\right) $ where each $f_{k}\in L_{s,0}^{2}\left(
\nu ^{2}\right) $ is a piecewise constant function with support in a
product set $S_{k}\times S_{k}$ with $\nu ^{2}$-finite measure. The
following isometric relation is well-known: $\forall f_{1},f_{2}\in
L_{s,0}^{2}\left( \nu
^{2}\right) $%
\begin{multline}
\mathds{E}\left[ I_{2}^{\tilde{N}^{c}}\left( f_{1}\right) \times I_{2}^{%
\tilde{N}^{c}}\left( f_{2}\right) \right]  \label{isometry} \\
=2\int_{\mathds{R}^{+}\times \mathbb{X}}\int_{\mathds{R}^{+}\times \mathbb{X}%
}f_{1}\left( s,x;t,y\right) f_{2}\left( s,x;t,y\right) \nu \left( \mathrm{d}%
s,\mathrm{d}x\right) \nu \left( \mathrm{d}t,\mathrm{d}y\right)
\text{.}
\end{multline}%
When $f\in L_{s}^{2}\left( \nu ^{2}\right) $ (hence $f$ does not
necessarily
vanish on diagonals), we set $I_{2}^{\tilde{N}^{c}}\left( f\right) =I_{2}^{%
\tilde{N}^{c}}\left( f \, \mathbb{I}_{D_{0}^{2}}\right) $, and we
observe that the isometry property (\ref{isometry}) still holds.
Indeed, $\nu $ is
non-atomic, and therefore $\nu ^{2}$ does not charge diagonals (even though $%
\tilde{N}^{c}\left( \mathrm{d}s,\mathrm{d}x\right)
\tilde{N}^{c}\left( \mathrm{d}t,\mathrm{d}y\right) $ does). We also
recall the \textsl{product
formula}%
\begin{multline}
\tilde{N}^{c}\left( g\right) \tilde{N}^{c}\left( h\right)
\label{eq : product} \\
=\left( g,h\right) _{L^{2}\left( \nu \right)
}+\int_{\mathds{R}^{+}\times
\mathbb{X}}g\left( s,x\right) h\left( s,x\right) \tilde{N}^{c}\left( \mathrm{%
d}s,\mathrm{d}x\right) +I_{2}^{\tilde{N}^{c}}\left( \widetilde{h\otimes g}%
\right) \text{,}
\end{multline}%
where $h\otimes g\left( s,x;t,y\right) =h\left( s,x\right) g\left(
t,y\right) \in L^{2}\left( \nu ^{2}\right) $ and ( $\widetilde{}$ )
stands for a symmetrization, which holds for every $f,g\in
L^{2}\left( \nu \right) $ such that $g\left( s,x\right) h\left(
s,x\right) \in L^{2}\left( \nu \right) $.

Finally, we state the main results proved in Peccati and Taqqu
(2006b). We consider a sequence of double integrals
\begin{equation}
F_{n}=I_{2}^{\tilde{N}^{c}}\left( f_{n}\right) \text{, \ \ }n\geq
1\text{,} \label{seq}
\end{equation}%
where $f_{n}\in L_{s,0}^{2}\left( \nu ^{2}\right) $. We will suppose
that
the following technical assumptions are satisfied: the sequence $f_{n}$, $%
n\geq 1$, in (\ref{seq}) is such that, for every $n\geq 1$,%
\begin{align}
& \left\Vert f_{n}\right\Vert _{L^{2}\left( \nu ^{2}\right)
}>0\text{ \ \ and \ }f_{n}\star _{2}^{1}f_{n}\in L^{2}\left( \nu
\right) ,  \tag{N1}
\label{N-1} \\
& \left\{ \int_{\mathds{R}^{+}\times \mathbb{X}}f_{n}\left(
s,y;\cdot
\right) ^{4}\nu \left( \mathrm{d}s,\mathrm{d}y\right) \right\} ^{\frac{1}{2}%
}\in L^{1}\left( \nu \right) ,  \tag{N2}  \label{N-2}
\end{align}%
where we use the notation introduced in (\ref{eq : f10})-(\ref{eq :
f21}), and moreover, as $n\rightarrow +\infty $,
\begin{equation}
\int_{\mathds{R}^{+}\times \mathbb{X}}\int_{\mathds{R}^{+}\times \mathbb{X}%
}f_{n}\left( s,y;t,x\right) ^{4}\nu \left(
\mathrm{d}s,\mathrm{d}y\right) \nu \left(
\mathrm{d}t,\mathrm{d}x\right) \rightarrow 0.  \tag{N3} \label{N-3}
\end{equation}
Note that (\ref{N-3}) implies, in particular, that $f_{n}\in
L^{4}\left( \nu ^{2}\right) $ for every $n$. See Peccati and Taqqu
(2006b) for a discussion of the role of (\ref{N-1})-(\ref{N-3}). In
the subsequent sections, we will see how such assumptions restrict
the set of the random hazard rates that can be studied by our
techniques. The next result is a CLT involving sequences of double
integrals.

\begin{theorem}[Peccati and Taqqu, 2006b, Th.7]
\label{T : PoissCLT}Define the sequence
$F_{n}=I_{2}^{\tilde{N}^{c}}(f_{n})$
and $f_{n}\in L_{s,0}^{2}(\nu ^{2})$, $n\geq 1$, as in \textnormal{%
\eqref{seq}}, and suppose \textnormal{\eqref{N-1}--\eqref{N-3}} hold. Then, $%
f_{n}\star _{1}^{0}f_{n}\in L^{2}(\nu ^{3})$ for every $n\geq 1$,
and moreover:

\begin{enumerate}
\item[$1.$] if
\begin{eqnarray}
\left\Vert f_{n}\right\Vert _{L^{2}(\nu ^{2})}^{-2}\times \left(
f_{n}\star _{1}^{1}f_{n}\right) &\rightarrow &0\text{ in }L^{2}(\nu
^{2})\text{ \ and \
}  \label{G*} \\
\left\Vert f_{n}\right\Vert _{L^{2}(\nu ^{2})}^{-2}\times \left(
f_{n}\star _{2}^{1}f_{n}\right) &\rightarrow &0\text{ in }L^{2}(\nu
)  \nonumber
\end{eqnarray}%
then
\begin{equation}
2^{-1/2}\left\Vert f_{n}\right\Vert _{L^{2}(\nu ^{2})}^{-1}\times F_{n}%
\overset{\mathrm{law}}{\longrightarrow}X\text{,}  \label{GG}
\end{equation}%
where $X\sim \mathscr{N}\left( 0,1\right) $ is a standard Gaussian
random variable;

\item[$2.$] if $F_{n}\in L^{4}\left( \mathds{P}\right) $ for every $n$, then
a sufficient condition to have \textnormal{\eqref{G*}} is that%
\begin{equation}
\left( 2\left\Vert f_{n}\right\Vert _{L^{2}(\nu ^{2})}^{2}\right) ^{-2}%
\mathds{E}\left( F_{n}^{4}\right) \rightarrow 3;  \label{GGG}
\end{equation}

\item[$3.$] if the sequence $\left\{ \left( 2\left\Vert f_{n}\right\Vert
_{L^{2}(\nu ^{2})}^{2}\right) ^{-2}F_{n}^{4}:n\geq 1\right\} $ is
uniformly integrable, then conditions \textnormal{\eqref{G*},
\eqref{GG}} and \textnormal{\eqref{GGG}} are equivalent.
\end{enumerate}
\end{theorem}

Theorem \ref{T : PoissCLT} is proved by using a decoupling
technique, known as the \textsl{principle of conditioning}, which
has been adapted to the framework of CRM by means of the general
theory of stable convergence developed in Peccati and Taqqu (2006a).
The next result gives sufficient conditions to have that the law of
a random vector, composed of a single and of a double integral,
converges weakly to a bivariate Gaussian law. The proof is
essentially based on an appropriate version of the \textsl{product
formulae }for multiple stochastic integrals, proved e.g. in
Surgailis (1984).

\begin{theorem}[Peccati and Taqqu, 2006b, Th. 8]
\label{T : jointCV} \ \newline \indent \textnormal{(A)} Consider a
sequence
\begin{equation}
G_{n}=\tilde{N}^{c}\left( g_{n}\right) ,\text{ \ \ }n\geq 1,
\nonumber
\end{equation}%
where $g_{n}\in L^{2}(\nu )\cap L^{3}(\nu )$ and $\left\Vert
g_{n}\right\Vert _{L^{2}\left( \nu \right) }>0$, and suppose that, as $%
n\rightarrow +\infty $,
\begin{equation}
\left\Vert g_{n}\right\Vert _{L^{2}\left( \nu \right) }^{-3}\int_{\mathds{R}%
^{+}\times \mathbb{X}}\left\vert g_{n}\left( s,y\right) \right\vert
^{3}\nu \left( \mathrm{d}s,\mathrm{d}y\right) \rightarrow 0\text{.}
\label{CVgn}
\end{equation}%
Then, $\left\Vert g_{n}\right\Vert _{L^{2}(\nu )}^{-1}\times G_{n}\overset{%
\mathrm{law}}{\longrightarrow}X$, where $X\sim \mathscr{N}\left(
0,1\right) $ is a centered standard Gaussian random variable.

\textnormal{(B)} Consider a sequence $F_{n}=I_{2}^{\tilde{N}^{c}}(f_{n})$, $%
n\geq 1$, with $f_{n}\in L_{s,0}^{2}(\nu ^{2})$ as in \textnormal{\eqref{seq}%
}, and a sequence $G_{n}=\tilde{N}^{c}\left( g_{n}\right) $, $n\geq
1$, as at Point \textnormal{(A)}. Suppose moreover that

\begin{description}
\item[\textnormal{(i)}] The sequence $\left( f_{n}\right) $ verifies
assumptions \textnormal{\eqref{N-1}--\eqref{N-3}}, and satisfies
condition \textnormal{\eqref{G*}};

\item[\textnormal{(ii)}] The sequence $\left( g_{n}\right) $ satisfies
\textnormal{\eqref{CVgn}}.
\end{description}

Then, as $n\rightarrow +\infty $,
\begin{equation}
\left( 2^{-1/2}\left\Vert f_{n}\right\Vert _{L^{2}(\nu
^{2})}^{-1}\times F_{n},\left\Vert g_{n}\right\Vert _{L^{2}(\nu
)}^{-1}\times G_{n}\right)
\overset{\mathrm{law}}{\longrightarrow}\left( X,X^{\prime }\right)
\text{,} \label{Jclt}
\end{equation}%
where $X,X^{\prime }\sim \mathscr{N}\left( 0,1\right) $ are two
independent, centered standard Gaussian random variables.
\end{theorem}

Part B of Theorem \ref{T : jointCV} implies in particular that,
whenever conditions (\ref{G*}) and (\ref{Jclt}) are met, the
(componentwise)
convergence of $\left\Vert f_{n}\right\Vert ^{-1}\times F_{n}$ and $%
\left\Vert g_{n}\right\Vert ^{-1}\times G_{n}$, towards a Gaussian
distribution, \textsl{implies necessarily }the joint convergence of
the vector $\left( \left\Vert f_{n}\right\Vert ^{-1}F_{n},\left\Vert
g_{n}\right\Vert ^{-1}G_{n}\right) $. This conclusion echoes results
already established in the framework of Gaussian CRM (see Peccati
and Tudor (2005)).

\smallskip

Now consider the positive kernel $k$, which defines $\tilde{h}$ via (\ref%
{eq:mixture_hazard}), and suppose (here and for the remainder of the
Section) that $k$ satisfies assumption (\ref{eq:H2}). In the next
two Lemmas we collect some straightforward facts which will be used
throughout the sequel.

\begin{lemma}
\label{L : Law}The two processes $\tilde{h}\left( t\right) $, $t\geq
0$, and
\begin{equation}
\tilde{h}_{\ast }\left( t\right) :=\tilde{N}^{c}\left( (\cdot
)k\left( t,\cdot \right) \right) +\int_{\mathds{R}^{+}\times
\mathbb{X}}sk\left( t,x\right) \nu \left(
\mathrm{d}s,\mathrm{d}x\right) ,\text{ \ \ }t\geq 0, \nonumber
\end{equation}%
where
\begin{equation}
\tilde{N}^{c}\left( (\cdot )k\left( t,\cdot \right) \right) :=\int_{%
\mathds{R}^{+}\times \mathbb{X}}sk\left( t,x\right)
\tilde{N}^{c}\left( \mathrm{d}s,\mathrm{d}x\right) ,  \label{nc1}
\end{equation}%
have the same law.
\end{lemma}

\textsc{Proof. }Use (\ref{eq:Laplace}) and (\ref{PLK}) to compute
the two
transforms%
\[
\mathds{E}\left[ \mathrm{e}^{i\sum_{j=1}^{n}\lambda
_{j}\tilde{h}\left(
t_{j}\right) }\right] \text{ \ and \ }\mathds{E}\left[ \mathrm{e}%
^{i\sum_{j=1}^{n}\lambda _{j}\tilde{h}_{\ast }\left( t_{j}\right)
}\right] ,
\]%
for every $n\geq 1$, every $\left( \lambda _{1},...,\lambda _{n}\right) \in %
\mathds{R}^{n}$ and every $t_{1},...,t_{n}\geq 0$. $\square $


\begin{lemma}
\label{L : Fubini}For every $T>0$,%
\begin{eqnarray}
\int_{0}^{T}\int_{\mathds{R}^{+}\times \mathbb{X}}sk\left( t,x\right) \tilde{%
N}^{c}\left( \mathrm{d}s,\mathrm{d}x\right) \mathrm{d}t &=&\tilde{N}%
^{c}\left( k_{T}^{\left( 0\right) }\right) ,  \label{F1} \\
\frac{1}{T}\int_{0}^{T}\int_{\mathds{R}^{+}\times
\mathbb{X}}s^{2}k\left(
t,x\right) ^{2}\tilde{N}^{c}\left( \mathrm{d}s,\mathrm{d}x\right) \mathrm{d}%
t &=&\tilde{N}^{c}\left( k_{T}^{\left( 2\right) }\right) ,
\label{F1b}
\end{eqnarray}%
where $k_{T}^{\left( 0\right) }$ and $k_{T}^{\left( 2\right) }$ are
given,
respectively, by \textnormal{\eqref{kaT1}} and \textnormal{\eqref{k2}}. If $%
k_{T}^{\left( 3\right) }\in L^{2}\left( \nu \right) $ $\cap
L^{1}\left( \nu
\right) $%
\begin{equation}
\frac{1}{T}\int_{0}^{T}\tilde{N}^{c}\left( (\cdot )k\left( t,\cdot
\right) \right) \left( \int_{\mathds{R}^{+}\times
\mathbb{X}}sk\left( t,x\right) \nu
\left( \mathrm{d}s,\mathrm{d}x\right) \right) \mathrm{d}t=\tilde{N}%
^{c}\left( k_{T}^{\left( 3\right) }\right) .  \label{F1c}
\end{equation}%
Analogously, for every $T>0$,%
\begin{equation}
\frac{1}{T}\int_{0}^{T}I_{2}^{\tilde{N}^{c}}\left( \left[ \left(
\cdot \right) k\left( t,\cdot \right) \right] \otimes \left[ \left(
\cdot \right)
k\left( t,\cdot \right) \right] \right) \mathrm{d}t=I_{2}^{\tilde{N}%
^{c}}\left( k_{T}^{\left( 1\right) }\right) ,  \label{F2}
\end{equation}%
where $\left[ \left( \cdot \right) k\left( t,\cdot \right) \right] \otimes %
\left[ \left( \cdot \right) k\left( t,\cdot \right) \right] \left(
u,x;v,y\right) :=uvk\left( t,x\right) k\left( t,y\right) $, and $%
k_{T}^{\left( 1\right) }$ is defined according to
\textnormal{\eqref{k1}}.
\end{lemma}

The proof of Lemma \ref{L : Fubini} is trivial when the map $\left(
t,x\right) \mapsto k\left( t,x\right) $ is piecewise constant:
indeed, in this case (\ref{F1}), (\ref{F1b}), (\ref{F1c}) and
(\ref{F2}) follow immediately from the application of a standard
Fubini theorem. The general statement is obtained by a density
argument; we omit the details here (one can e.g. mimic the proof of
Lemma 13 in Peccati, 2001).

Finally note that, given two sequences of random variables $\left\{
A_{n}\right\} $ and $\left\{ B_{n}\right\}$ such that $A_{n}-B_{n}%
\rightarrow 0$ in probability, we will sometimes write
\begin{equation}
A_{n}\overset{\mathds{P}}{\approx }B_{n}\text{.}  \nonumber
\end{equation}

\subsection{Proof of Theorem \protect\ref{T : cumHRasy}}

Use Lemma \ref{L : Law} and relations (\ref{nc1}) and (\ref{F1}) to
write
\begin{align*}
\tilde{H}(T)&\overset{\mathrm{law}}{=}\int_{0}^{T}\tilde{h}_{\ast
}\left(
t\right) \mathrm{d}t \\
&=\int_{0}^{T}\tilde{N}^{c}\left( (\cdot )k\left( t,\cdot \right)
\right) \mathrm{d}t+\int_{0}^{T}\int_{\mathds{R}^{+}\times
\mathbb{X}}sk\left(
t,x\right) \nu \left( \mathrm{d}s,\mathrm{d}x\right) \mathrm{d}t \\
&=\int_{0}^{T}\int_{\mathds{R}^{+}\times \mathbb{X}}sk\left(
t,x\right)
\tilde{N}^{c}\left( \mathrm{d}s,\mathrm{d}x\right) \mathrm{d}%
t+\int_{0}^{T}\int_{\mathds{R}^{+}\times \mathbb{X}}sk\left(
t,x\right) \nu
\left( \mathrm{d}s,\mathrm{d}x\right) \mathrm{d}t \\
&=\tilde{N}^{c}\left( k_{T}^{\left( 0\right) }\right) +\int_{0}^{T}\int_{%
\mathds{R}^{+}\times \mathbb{X}}sk\left( t,x\right) \nu \left( \mathrm{d}s,%
\mathrm{d}x\right) \mathrm{d}t,
\end{align*}%
which yields, via the relation $\mathds{E}(\tilde{H}(T))=\int_{0}^{T}\int_{%
\mathds{R}^{+}\times \mathbb{X}}sk\left( t,x\right) \nu \left( \mathrm{d}s,%
\mathrm{d}x\right) \mathrm{d}t$,
\[
C_{0}\left( k,T\right) \times \left[ \tilde{H}(T)-\mathds{E}(\tilde{H}(T))%
\right] \overset{\mathrm{law}}{=}\tilde{N}^{c}\left( C_{0}\left(
k,T\right) \times k_{T}^{\left( 0\right) }\right) \text{.}
\]%
Since the isometry property (\ref{eq : isoPoiss}) and the assumption (\ref%
{a-}) yield%
\[
\mathds{E}\lbrack \tilde{N}^{c}(C_{0}\left( k,T\right) \times
k_{T}^{\left( 0\right) })^{2}] =C_{0}^{2}\left( k,T\right)
\int_{\mathds{R}^{+}\times \mathbb{X}}\!\left[ k_{T}^{\left(
0\right) }\left( s,x\right) \right] ^{2}\!\nu \left(
\mathrm{d}s,\mathrm{d}x\right) \rightarrow \sigma _{0}^{2}\left(
k\right) ,
\]%
we deduce from Part A of Theorem \ref{T : jointCV} (in the case $%
g_{n}=(C_{0}\left( k,T_{n}\right) /\sigma _{0}\left( k\right)
)\times k_{T_{n}}^{\left( 0\right) }$, where $T_{n}$ is any positive
sequence diverging to infinity) that, since (\ref{aa}) holds, the
CLT (\ref{aaa}) must also take place. \hfill \qed

\subsection{Proof of Theorem \protect\ref{T : path2ndM}}

Use Lemma \ref{L : Law} to write (we adopt once again the notation (\ref{nc1}%
))%
\begin{align*}
\frac{1}{T}\!\int_{0}^{T}\tilde{h}(t)^{2}\mathrm{d}t& \overset{\mathrm{law}}{%
=}\frac{1}{T}\!\int_{0}^{T}\!\!\tilde{N}^{c}\left( (\cdot )k\left(
t,\cdot
\right) \right) ^{2}\mathrm{d}t+\frac{1}{T}\!\int_{0}^{T}\!\!\left( \int_{%
\mathds{R}^{+}\!\times \mathbb{X}}\!\!\!\!sk\left( t,x\right) \nu
\left(
\mathrm{d}s,\mathrm{d}x\right) \right) ^{2}\!\!\mathrm{d}t \\
& \qquad \qquad \ \ +\frac{2}{T}\!\int_{0}^{T}\tilde{N}^{c}\left(
(\cdot
)k\left( t,\cdot \right) \right) \left( \int_{\mathds{R}^{+}\times \mathbb{X}%
}\!\!\!sk\left( t,x\right) \nu \left( \mathrm{d}s,\mathrm{d}x\right)
\right) \mathrm{d}t.
\end{align*}%
Now recall that, thanks to (\ref{F1c}),%
\[
\frac{2}{T}\int_{0}^{T}\tilde{N}^{c}\left( (\cdot )k\left( t,\cdot
\right) \right) \left( \int_{\mathds{R}^{+}\times
\mathbb{X}}sk\left( t,x\right) \nu
\left( \mathrm{d}s,\mathrm{d}x\right) \right) \mathrm{d}t=\tilde{N}%
^{c}\left( 2k_{T}^{\left( 3\right) }\right) ,
\]%
so that, by using (\ref{eq : semplVar}),
\begin{eqnarray}
&&C_{1}\left( k,T\right) \times \left\{ \frac{1}{T}\int_{0}^{T}\tilde{h}%
(t)^{2}\mathrm{d}t-\frac{1}{T}\int_{0}^{T}\mathds{E}\lbrack \tilde{h}(t)^{2}]%
\mathrm{d}t\right\}   \nonumber \\
&&\overset{\mathrm{law}}{=}C_{1}\left( k,T\right) \times \left\{ \frac{1}{T}%
\int_{0}^{T}\tilde{N}^{c}\left( (\cdot )k\left( t,\cdot \right) \right) ^{2}%
\mathrm{d}t^{^{^{^{^{^{{}}}}}}}\right.   \nonumber \\
&&\text{ \ \ \ \ \ \ \ \ \ \ \ \ \ \ \ \ \ \ \ \ \ \ \ }\left. +\text{ }%
\tilde{N}^{c}\left( 2k_{T}^{\left( 3\right) }\right) -\frac{1}{T}%
\int_{0}^{T}\int_{\mathds{R}^{+}\times \mathbb{X}}\!\!s^{2}k\left(
t,x\right) ^{2}\nu \left( \mathrm{d}s,\mathrm{d}x\right) \right\} \mathrm{d}%
t.  \label{quad1}
\end{eqnarray}%
By applying the product formula (\ref{eq : product}) in the case
$g\left( s,x\right) =h\left( s,x\right) =sk\left( t,x\right) $, for
every $t\geq 0$ we obtain
\begin{multline*}
\tilde{N}^{c}\left( (\cdot )k\left( t,\cdot \right) \right) ^{2}=\int_{%
\mathds{R}^{+}\times \mathbb{X}}s^{2}k\left( t,x\right) ^{2}\nu
\left(
\mathrm{d}s,\mathrm{d}x\right)  \\
+\int_{\mathds{R}^{+}\times \mathbb{X}}s^{2}k\left( t,x\right) ^{2}\tilde{N}%
^{c}\left( \mathrm{d}s,\mathrm{d}x\right)
+I_{2}^{\tilde{N}^{c}}\left( \left[ \left( \cdot \right) k\left(
t,\cdot \right) \right] \otimes \left[ \left( \cdot \right) k\left(
t,\cdot \right) \right] \right) \text{,}
\end{multline*}%
from which we deduce that, thanks to formulae (\ref{F1b}) and
(\ref{F2}), the expression in (\ref{quad1}) is indeed equal to
\begin{equation}
C_{1}\left( k,T\right) \times \left\{ \tilde{N}^{c}\left(
k_{T}^{\left( 2\right) }+2k_{T}^{\left( 3\right) }\right)
+I_{2}^{\tilde{N}^{c}}\left( k_{T}^{\left( 1\right) }\right)
\right\} ,  \nonumber
\end{equation}%
for every $T>0$. It follows that Theorem \ref{T : path2ndM} is
proved, once it is shown that
\[
\left( \tilde{N}^{c}\left( C_{1}\left( k,T\right) \times
(k_{T}^{\left( 2\right) }+2k_{T}^{\left( 3\right) })\right)
,I_{2}^{\tilde{N}^{c}}\left( C_{1}\left( k,T\right) \times
k_{T}^{\left( 1\right) }\right) \right)
\overset{\mathrm{law}}{\longrightarrow }\left( X,X^{\prime }\right)
\]%
where $X$ and $X^{\prime }$ are independent and such that $X\sim \mathscr{N}%
\left( 0,\sigma _{2}^{2}\left( k\right) \right) $ and $X^{\prime }\sim %
\mathscr{N}\left( 0,\sigma _{1}^{2}\left( k\right) \right) $. To
this end, we apply Part B of Theorem \ref{T : jointCV}: according to
such a result, it is sufficient to check that, for every positive
sequence $T_{n}\rightarrow +\infty $, the two sequences
\[
g_{n}=\frac{C_{1}\left( k,T_{n}\right) }{\sigma _{2}\left( k\right) }%
(k_{T_{n}}^{\left( 2\right) }+2k_{T_{n}}^{\left( 3\right) })\text{ \
\ and \
\ }f_{n}=\frac{C_{1}\left( k,T_{n}\right) }{\sigma _{1}\left( k\right) }%
k_{T_{n}}^{\left( 1\right) }\text{,\ \ }n\geq 1\text{,}
\]%
satisfy, respectively, condition (\ref{CVgn}) and conditions (\ref{N-1})-(%
\ref{N-3}) and (\ref{G*}). It is immediately seen that Assumptions 5
and 6 in the statement imply (\ref{CVgn}), and we are therefore left
with the sequence $\left\{ f_{n}\right\} $. Conditions (\ref{N-1})
and (\ref{N-2}) can be checked by standard iterations of the Jensen
and Cauchy-Schwarz inequalities (see e.g. Section 5.1 in Peccati and
Taqqu (2006b) for several analogous computations). Finally,
(\ref{N-3}) is given by Assumption 2 in the statement, whereas
Assumptions 3 and 4 give, respectively, the first and the second
line in (\ref{G*}). This concludes the proof of Theorem \ref{T :
path2ndM}. \hfill \qed

\subsection{Proof of Theorem \protect\ref{P : path variance}}

Write first
\begin{equation}
\frac{1}{T}\int_{0}^{T}[\tilde{h}(t)-\frac{1}{T}\tilde{H}(T)]^{2}\mathrm{d}t=%
\frac{1}{T}\int_{0}^{T}\tilde{h}(t)^{2}\mathrm{d}t-\left( \frac{1}{T}\tilde{H%
}(T)\right) ^{2},  \label{x}
\end{equation}%
and observe that
\begin{eqnarray}
C_{1}\left( k,T\right) \left( \frac{1}{T}\tilde{H}(T)\right) ^{2} &=&\frac{%
C_{1}\left( k,T\right) }{T^{2}C_{0}\left( k,T\right) ^{2}}\left\{
C_{0}\left( k,T\right) \left[
\tilde{H}(T)-\mathds{E}(\tilde{H}(T))\right]
\right\} ^{2}  \nonumber \\
&& \ +\, \frac{C_{1}\left( k,T\right) }{T^{2}}\,
\mathds{E}(\tilde{H}(T))^{2}
\label{xx} \\
&&\ +\,2\,\frac{C_{1}\left( k,T\right) }{T^{2}}\,\mathds{E}(\tilde{H}(T))%
\left[ \tilde{H}(T)-\mathds{E}(\tilde{H}(T))\right] .  \nonumber
\end{eqnarray}%
From Assumption 1 in the statement, and since (\ref{a-}) and
(\ref{aa}) are in order, we deduce
\begin{equation}
\frac{C_{1}\left( k,T\right) }{T^{2}C_{0}\left( k,T\right)
^{2}}\left\{ C_{0}\left( k,T\right) \left[
\tilde{H}(T)-\mathds{E}(\tilde{H}(T))\right] \right\}
^{2}\overset{\mathds{P}}{\rightarrow }0.  \label{xxx}
\end{equation}%
Moreover, Assumption 2 in the statement yields that, as
$T\rightarrow +\infty $,
\begin{multline}
\frac{2C_{1}\left( k,T\right)}{T^2} \mathds{E}(\tilde{H}(T))\left[ \tilde{H}%
(T)-\mathds{E}(\tilde{H}(T))\right] \\
\overset{\mathds{P}}{\approx }\delta \left( k\right) C_{0}\left( k,T\right) %
\left[ \tilde{H}(T)-\mathds{E}(\tilde{H}(T))\right]  \label{xxxx}
\end{multline}%
In view of Lemma \ref{L : Law}, and by reasoning as in the proof of Theorem %
\ref{T : cumHRasy} and Theorem \ref{T : path2ndM}, we infer from relations (%
\ref{x})-(\ref{xxxx}) that
\begin{align*}
& C_{1}\left( k,T\right) \times \left\{ \frac{1}{T}\int_{0}^{T}[\tilde{h}(t)-%
\frac{1}{T}\tilde{H}(T)]^{2}\mathrm{d}t-\frac{1}{T}\int_{0}^{T}\mathds{E}%
\lbrack \tilde{h}(t)^{2}]\mathrm{d}t+\frac{\mathds{E}\lbrack \tilde{H}%
(T))]^{2}}{T^{2}}\right\} \\
& \overset{\mathrm{law}}{=}\!\!\tilde{N}^{c}\!\left( C_{1}\left(
k,T\right)
\left( \!k_{T}^{\left( 2\right) }+2k_{T}^{\left( 3\right) }\right) \!-\!%
\frac{2C_{1}\left( k,T\right)
}{T^{2}}\mathds{E}(\tilde{H}(T))k_{T}^{\left( 0\right) }\!\right)
\!\!+\!\!I_{2}^{\tilde{N}^{c}}\left( \!C_{1}\left(
k,T\right) k_{T}^{\left( 1\right) }\!\right) \\
& \overset{\mathds{P}}{\approx }\tilde{N}^{c}\left( C_{1}\left(
k,T\right) \left( k_{T}^{\left( 2\right) }+2k_{T}^{\left( 3\right)
}\right) -\delta \left( k\right) C_{0}\left( k,T\right)
k_{T}^{\left( 0\right) }\right) +I_{2}^{\tilde{N}^{c}}\left(
C_{1}\left( k,T\right) k_{T}^{\left( 1\right) }\right) \text{.}
\end{align*}%
The conclusion is deduced from Assumption 3 in the statement, by
applying Theorem \ref{T : jointCV} in the case
\begin{eqnarray*}
g_{n} &=&\frac{C_{1}\left( k,T_{n}\right) \left( k_{T_{n}}^{\left(
2\right) }+2k_{T_{n}}^{\left( 3\right) }\right) -\delta \left(
k\right) C_{0}\left(
k,T_{n}\right) k_{T_{n}}^{\left( 0\right) }}{\sigma _{3}\left( k\right) } \\
f_{n} &=&\frac{C_{1}\left( k,T_{n}\right) }{\sigma _{1}\left( k\right) }%
k_{T_{n}}^{\left( 1\right) }\text{, \ \ }n\geq 1\text{,}
\end{eqnarray*}%
where $T_{n}\rightarrow +\infty .$ \hfill \qed

\subsection{Proof of Theorem \protect\ref{T : Comparison}}

To prove Part (A), observe that the assumptions imply the existence
of two
constants $0<D_{1}<D_{2}<+\infty $, such that, for $T$ sufficiently large,%
\[
D_{1}<C_{0}^{2}\left( k,T\right) \times \int_{\mathds{R}^{+}\times \mathbb{X}%
}\left[ k_{T}^{\left( 0\right) }\left( s,x\right) \right] ^{2}\nu
\left( \mathrm{d}s,\mathrm{d}x\right) <D_{2}\text{.}
\]%
Standard arguments yield therefore that, for every sequence $%
T_{n}\rightarrow +\infty $, there exists a subsequence $T_{n^{\prime
}}$
such that, as $n^{\prime }\rightarrow +\infty $,%
\[
C_{0}^{2}\left( k,T_{n^{\prime }}\right) \times
\int_{\mathds{R}^{+}\times
\mathbb{X}}\left[ k_{T_{n^{\prime }}}^{\left( 0\right) }\left( s,x\right) %
\right] ^{2}\nu \left( \mathrm{d}s,\mathrm{d}x\right) \rightarrow
\sigma ^{2}\left( k\right) >0\text{,}
\]%
where $\sigma ^{2}\left( k\right) $ is some well chosen positive
constant. Moreover,
\begin{align*}
C_{0}^{3}\left( k,T_{n^{\prime }}\right) \times &
\int_{\mathds{R}^{+}\times
\mathbb{X}}\left[ k_{T_{n^{\prime }}}^{\left( 0\right) }\left( s,x\right) %
\right] ^{3}\nu \left( \mathrm{d}s,\mathrm{d}x\right)  \\
& \leq C_{0}^{3}\left( k,T_{n^{\prime }}\right)
\int_{\mathds{R}^{+}}\left[ \overline{\overline{k}}_{T_{n^{\prime
}}}^{\left( 0\right) }\left(
s,x\right) \right] ^{3}\nu \left( \mathrm{d}s,\mathrm{d}x\right)  \\
& \sim C_{0}^{3}\left( \overline{\overline{k}},T_{n^{\prime }}\right) \int_{%
\mathds{R}^{+}}\left[ \overline{\overline{k}}_{T_{n^{\prime
}}}^{\left(
0\right) }\left( s,x\right) \right] ^{3}\nu \left( \mathrm{d}s,\mathrm{d}%
x\right) \rightarrow 0.
\end{align*}%
The proofs of Parts (B) and (C) are based on analogous computations,
and are omitted. \hfill \qed

\section{Conclusions and future work}

(I) Future research will focus on the generalization of our
asymptotic results to general multiplicative intensity models
(Aalen, 1978), which include a wide variety of popular models such
as Cox proportional hazards regression models, multiple decrement
models, birth and death processes and non--homogeneous Poisson
processes. To fix ideas consider the Cox
proportional hazards regression model, in which $Z_{i}$ is an $m$%
--dimensional vector of covariates recorded for the $i$--th individual and $%
\theta $ is a $m$--dimensional vector of unknown regression
coefficients. Then the proportional hazards model is specified in
terms of the hazard function relationship as
\begin{equation}
h_{i}(t)=h_{0}(t)\exp (\theta ^{\prime }\,Z_{i}),  \nonumber
\end{equation}%
where $h_{0}$ represents the so--called baseline hazard function. A
Bayesian treatment leads to considering $h_{0}$ and $\theta $ to be
random and,
hence, by choosing $\tilde{h}_{0}$ to be a mixture as in %
\eqref{eq:mixture_hazard} and $\pi $ to be a prior for
$\tilde{\theta}$, one obtains a semi--parametric random hazard rate
function for the $i$--th individual of the form
\begin{equation}
\tilde{h}_{i}(t)=\exp (\tilde{\theta}^{\prime }\,Z_{i})\int_{\mathbb{X}%
}k(t,x)\tilde{\mu}(\mathrm{d}x).  \label{eq:cox_hazard}
\end{equation}%
Bayesian analysis of the Cox model within this setup has been
pursued in Ibrahim, Chen and Mac Eachern (1999), James (2003),
Ishwaran and James (2004), Nieto--Barajas and Walker (2005). Since
\eqref{eq:mixture_hazard} still represents the basic building block
of \eqref{eq:cox_hazard} and, indeed, also of other multiplicative
intensity models, we aim at extending our results to random objects
such as \eqref{eq:cox_hazard} and expect to obtain CLTs for which
the limiting random variable is a suitable mixture of Gaussian
distributions.

(II) The techniques exploited in Section \ref{S : Proofs}, for
deriving the main results of this paper, can be further generalized.
As already mentioned, they are indeed based on a very general
decoupling criterion, known as the \textsl{principle of
conditioning}. As shown in Peccati and Taqqu (2006a,b), this
principle can be applied to a wide class of stochastic integrals
with respect to completely random measures, including multiple
Wiener-It\^{o} integrals of any order $n>2$. In particular, we
expect that the results of the present paper can be suitably
extended to accommodate the asymptotic analysis of non-linear and
non-quadratic functionals, such as e.g. path-moments of order
greater than two. Note that results of this type\ are already
available in the Gaussian case. See, e.g., Peccati and Tudor (2005).

\smallskip


\begin{center}
\textbf{REFERENCES}
\end{center}

{\footnotesize \noindent \textsc{Aalen, O.} (1978). Nonparametric
Inference for a Family of Counting Processes. \textit{Ann. Statist.}
\textbf{6}, 701--726. }

{\footnotesize \noindent \textsc{Brix, A.} (1999). Generalized gamma
measures and shot-noise Cox processes. \textit{Adv. Appl. Prob.} \textbf{31}%
, 929--953. }

{\footnotesize \noindent \textsc{Cifarelli, D.M.} and
\textsc{Regazzini, E.} (1990). Distribution functions of means of a
Dirichlet process. \textit{Ann. Statist.}, \textbf{18}, 429--442. }

{\footnotesize \noindent \textsc{Cifarelli, D.M.} and
\textsc{Melilli, E.} (2000). Some new results for Dirichlet priors.
\textit{Ann. Statist.} \textbf{28}, 1390--1413. }

{\footnotesize \noindent \textsc{De Blasi, P. \textrm{and} Hjort, N.
L.} (2006). Bayesian survival analysis in proportional hazard models
with logistic relative risk. \textit{Scand. J. Statist.}, to
appear.}

{\footnotesize \noindent \textsc{Daley, D. \textrm{and} Vere-Jones,
D. J.} (1988). \textsl{An introduction to the theory of point
processes}. Springer, New York.}

{\footnotesize \noindent \textsc{Doksum, K.} (1974). Tailfree and
neutral
random probabilities and their posterior distributions. \textit{Ann. Probab.}%
, \textbf{2}, 183--201. }

{\footnotesize \noindent \textsc{Dykstra, R.L. \textrm{and} Laud,
P.} (1981). A Bayesian nonparametric approach to reliability.
\textit{Ann. Statist.} \textbf{9} 356---367. }

{\footnotesize \noindent \textsc{Epifani, I., Lijoi, A. \textrm{and}
Pr\"unster, I.} (2003). Exponential functionals and means of
neutral-to-the-right priors. \textit{Biometrika}, \textbf{90},
791--808. }

{\footnotesize \noindent \textsc{Ferguson, T.S.} (1974). Prior
distributions on spaces of probability measures. \textit{Ann.
Statist.}, \textbf{2}, 615--629. }

{\footnotesize \noindent \textsc{Ferguson, T.S. \textrm{and} Phadia,
E.G.}
(1979). Bayesian nonparametric estimation based on censored data. \textit{%
Ann. Statist.}, \textbf{7}, 163--186. }

{\footnotesize \noindent \textsc{Ghosh, J.K.~\textrm{and}
Ramamoorthi, R.V.} (2003). \textsl{Bayesian Nonparametrics.}
Springer, New York. }

{\footnotesize \noindent \textsc{Hjort, N. L.} (1990). Nonparametric
Bayes
estimators based on beta processes in models for life history data. \textit{%
Ann. Statist.} \textbf{18}, 1259--1294. }

{\footnotesize \noindent \textsc{Ho, M.-W.} (2006). A Bayes method
for a monotone hazard rate via S-paths. \textit{Ann. Statist.}
\textbf{34}, 820--836. }

{\footnotesize \noindent \textsc{Ibrahim, J.G., Chen, M.-H.
\textrm{and} MacEachern, S.N.} (1999). Bayesian variable selection
for proportional hazards models. \textit{Canad. J. Statist.}
\textbf{27}, 701--717. }

{\footnotesize \noindent \textsc{Ishwaran, H.} and \textsc{James, L.
F.} (2001). Gibbs sampling methods for stick-breaking priors.
\textit{J. Amer. Stat. Assoc.}, \textbf{96}, 161--173. }

{\footnotesize \noindent \textsc{Ishwaran, H.} and \textsc{James, L.
F.} (2004). Computational methods for multiplicative intensity
models using weighted gamma processes: Proportional hazards, marked
point processes, and panel count data. \textit{J. Amer. Stat.
Assoc.} \textbf{99}, 175--190. }

{\footnotesize \noindent \textsc{James L.F.} (2003). Bayesian
calculus for gamma processes with applications to semiparametric
intensity models. \textit{Sankhya} \textbf{65} (2003), 179--206. }

{\footnotesize \noindent \textsc{James L.F.} (2005). Bayesian
Poisson process partition calculus with an application to Bayesian
L\'{e}vy moving averages. \textit{Ann. Statist.} \textbf{33},
1771--1799. }

{\footnotesize \noindent \textsc{James L.F.} (2006). Poisson
calculus for spatial neutral to the right processes. \textit{Ann.
Statist.} \textbf{34}, 416--440. }

{\footnotesize \noindent \textsc{Kalbfleisch, J.D.} (1978).
Non-parametric
Bayesian analysis of survival time data. \textit{J. Roy. Statist. Soc. Ser. B%
} \textbf{40}, 214--221. }

{\footnotesize \noindent \textsc{Kim, Y.} (1999). Nonparametric
Bayesian estimators for counting processes. \textit{Ann. Statist.}
\textbf{27}, 562--588. }

{\footnotesize \noindent \textsc{Kim, Y. \textrm{and} Lee, J.}
(2003). Bayesian analysis of proportional hazard models.
\textit{Ann. Statist.} \textbf{31}, 493--511. }

{\footnotesize \noindent \textsc{Kingman, J.F.C.} (1967). Completely
random measures. \textit{Pacific J. Math.}, \textbf{21}, 59-78. }

{\footnotesize \noindent \textsc{Kingman, J.F.C.} (1993).
\textit{Poisson Processes}. Oxford University Press, Oxford. }

{\footnotesize \noindent \textsc{Kwapien, S. \textrm{and} Woyczynski, W.A. }%
(1992). \textit{Random Series and Stochastic Integrals: Single and Multiple}%
. Birkh\"{a}user, Basel.}

{\footnotesize \noindent \textsc{Laud, P., Smith, A.F.M.
\textrm{and} Damien, P.} (1996). Monte Carlo methods for
approximating a posterior hazard rate process. \textit{Stat.
Comput.} \textbf{6}, 77--83. }

{\footnotesize \noindent \textsc{Lijoi, A., Mena, R. \textrm{and}
Pr\"unster, I.} (2005). Hierarchical mixture modelling with
normalized inverse Gaussian priors. \textit{J. Amer. Stat. Assoc.}
\textbf{100}, 1278--1291. }

{\footnotesize \noindent \textsc{Lo, A. Y. \textrm{and} Weng,
C.--S.} (1989). On a class of Bayesian nonparametric estimates. II.
Hazard rate estimates. \textit{Ann. Inst. Statist. Math.}
\textbf{41}, 227---245. }

{\footnotesize \noindent \textsc{Nieto-Barajas, L.E. \textrm{and}
Walker, S.G.} (2004). Bayesian nonparametric survival analysis via
L\'evy driven Markov processes. \textit{Statist. Sinica}
\textbf{14}, 1127--1146. }

{\footnotesize \noindent \textsc{Nieto-Barajas, L.E. \textrm{and}
Walker,
S.G.} (2005). A semi-parametric Bayesian analysis of survival data based on L%
\'{e}vy-driven processes. \textit{Lifetime Data Anal.} \textbf{11},
529--543. }

{\footnotesize \noindent \textsc{Nualart, D. \textrm{and} Vives, J.}
(1990). Anticipative calculus for the Poisson process based on the
Fock space . \textit{S\'{e}minaire de Probabilit\'{e}s XXIV}, LNM
\textbf{1426}, 154--165, Springer, Berlin.}

{\footnotesize \noindent \textsc{Peccati, G.} (2001). On the
convergence of
multiple random integrals. \textit{Studia Scient. Math. Hungarica} \textbf{37%
}, 429-470.}

{\footnotesize \noindent \textsc{Peccati, G. \textrm{and} Taqqu,
M.S.} (2006a). Stable convergence of generalized $L^{2}\!$ integrals
and the principle of conditioning. Preprint available at
www.geocities .com/giovannipeccati }

{\footnotesize \noindent \textsc{Peccati, G. \textrm{and} Taqqu,
M.S.} (2006b). Central limit theorems for double Poisson integrals.
Preprint available at www.geocities.com/giovannipeccati }

{\footnotesize \noindent \textsc{Peccati, G. \textrm{and} Tudor,
C.A.} (2005). Gaussian limits for vector-valued multiple stochastic
integrals. \textit{S\'{e}minaire de Probabilit\'{e}s XXXVIII}, LNM
\textbf{1857}, 247--262, Springer, Berlin.}

{\footnotesize \noindent \textsc{Regazzini, E., Guglielmi, A.
\textrm{and} Di Nunno, G.} (2002). Theory and numerical analysis for
exact distribution of functionals of a Dirichlet process.
\textit{Ann. Statist.} \textbf{30}, 1376--1411. }

{\footnotesize \noindent \textsc{Regazzini, E., Lijoi, A. \textrm{and} Pr%
\"{u}nster, I.} (2003). Distributional results for means of random
measures with independent increments. \textit{Ann. Statist.}
\textbf{31}, 560--585. }

{\footnotesize \noindent \textsc{Rota, G.--C. \textrm{and}
Wallstrom, C.} (1997). Stochastic integrals: a combinatorial
approach. \textit{Ann. Prob.} \textbf{25(}3), 1257-1283.}

{\footnotesize \noindent \textsc{Sato, K.} (1999). \textit{L\'evy
Processes and Infinitely Divisible Distributions.} Cambridge
University Press, Cambridge.}

{\footnotesize \noindent \textsc{Surgailis, D.} (1984). On multiple
Poisson integrals and associated Markov semigroups. \textit{Prob.
Math. Statist.} \textbf{3}, 217-239.}

{\footnotesize \noindent \textsc{Surgailis, D.} (2000). Non-CLT's:
U-statistics, Multinomial Formula and Approximations of multiple Wiener-It%
\^{o} integrals. In: \textsl{Long Range Dependence.}, 129--142, Birk\"{a}%
user, Basel.}

{\footnotesize \noindent \textsc{Walker, S. \textrm{and} Damien, P.}
(1998). A full Bayesian non-parametric analysis involving a neutral
to the right process. \textit{Scand. J. Statist.}, \textbf{25},
669--680. }

{\footnotesize \noindent \textsc{Walker, S. \textrm{and} Muliere,
P.} (1997). Beta-Stacy processes and a generalization of the
P\'olya-urn scheme. \textit{Ann. Statist.}, \textbf{25}, 1762--1780.
}

\end{document}